\documentclass[11pt]{article}

\usepackage[utf8]{inputenc} 
\usepackage[T1]{fontenc}    
\usepackage{hyperref}       
\usepackage{url}            
\usepackage{booktabs}       
\usepackage{amsfonts}       
\usepackage{nicefrac}       
\usepackage{microtype}      

\usepackage{amsmath}
\usepackage{amssymb}
\usepackage{amsthm}
\usepackage{pifont}
\usepackage{algorithm}
\usepackage{algpseudocode}
\usepackage[capitalize]{cleveref}
\usepackage{breqn}
\usepackage{xspace}
\usepackage{forloop}
\usepackage{multirow}
\usepackage{bbm}
\usepackage[pdftex,dvipsnames,table]{xcolor}
\usepackage[colorinlistoftodos,prependcaption,textsize=small]{todonotes}
\usepackage{graphicx}
\usepackage{paralist}
\usepackage{subcaption}
\usepackage{dsfont}
\usepackage{natbib}
\usepackage{csquotes}

\newcommand{\sfw}{\texttt{1-SFW}\xspace}

\newtheorem{theorem}{Theorem}
\newtheorem{lemma}{Lemma}
\newtheorem{assump}{Assumption}
\newtheorem{remark}{Remark}

\crefname{assump}{Assumption}{Assumptions}

\newcommand{\tnabla}{\tilde{\nabla}}
\newcommand{\tF}{\tilde{F}}
\newcommand{\expect}{\mathbb{E}}
\newcommand{\constraint}{\mathcal{K}}

\newcommand{\bd}{\mathbf{d}}
\newcommand{\bx}{\mathbf{x}}

\newcommand{\bv}{\mathbf{v}}
\newcommand{\bu}{\mathbf{u}}
\newcommand{\by}{\mathbf{y}}
\newcommand{\bz}{\mathbf{z}}

\DeclareMathOperator*{\argmax}{arg\,max}
\DeclareMathOperator*{\argmin}{arg\,min}

\input{required.tex}

\begin{document}
\title{One Sample Stochastic Frank-Wolfe}

\author{\\
	Mingrui Zhang\\
	Yale University\\
	\texttt{mingrui.zhang@yale.edu}
	\and\\
	Zebang Shen \\ 
	University of Pennsylvania \\
	\texttt{zebang@seas.upenn.edu}
	\and \\
	Aryan Mokhtari\\
	The University of Texas at Austin \\
	\texttt{mokhtari@austin.utexas.edu}
	\and\\
	Hamed Hassani \\
	University of Pennsylvania\\
	\texttt{hassani@seas.upenn.edu}
	\and\\
	Amin Karbasi\\
	Yale University \\
	\texttt{amin.karbasi@yale.edu}\\
}

\maketitle

%

%

%
%
%

\begin{abstract}
\vspace{-2mm}
One of the beauties of the projected gradient descent method lies in its rather simple mechanism and yet stable  behavior with inexact, stochastic  gradients, which has led to its wide-spread use in many machine learning applications.  However, once we replace the projection operator with a simpler linear program, as is done in the Frank-Wolfe method, both simplicity and  stability take a serious hit. The aim of this paper is to bring them back without sacrificing the efficiency. 
In this paper, we propose the first one-sample stochastic Frank-Wolfe algorithm, called \sfw,  that avoids the need to carefully tune the batch size, step size, learning rate, and other complicated hyper parameters. In particular, \sfw achieves the optimal  convergence rate  of $\mathcal{O}(1/\epsilon^2)$  for reaching an  $\epsilon$-suboptimal solution in the stochastic convex setting, and  a $(1-1/e)-\epsilon$ approximate solution for a stochastic  monotone DR-submodular maximization problem. Moreover, in a general non-convex setting, \sfw finds an $\epsilon$-first-order stationary point after at most $\mathcal{O}(1/\epsilon^3)$ iterations, achieving the current  best known convergence rate. All of this is possible by designing a novel unbiased momentum estimator that governs the stability of the optimization process while using a single sample at each iteration.  
\vspace{-2mm}

\end{abstract}

\newpage

\section{Introduction}

Projection-free methods, also known as conditional gradient methods or 
Frank-Wolfe (FW) methods, have been widely used for solving constrained 
optimization problems 
\citep{frank1956algorithm,jaggi2013revisiting,lacoste2015global}.  Indeed, 
extending such methods to the stochastic setting is a challenging task as it is 
known that FW-type methods are highly sensitive to stochasticity 
in gradient computation \citep{DBLP:conf/icml/HazanK12}. To resolve this issue 
several stochastic variants of FW methods have been studied 
in the literature 
\citep{DBLP:conf/icml/HazanK12,hazan2016variance,reddi2016stochastic,lan2016conditional,braun2017lazifying,hassani2019stochastic,shen2019complexities,
	yurtsever2019conditional}. In all these stochastic methods, 
the basic idea is to provide an accurate estimate of the gradient by using some 
variance-reduction techniques that typically rely on large mini-batches of samples where the size grows with the 
number of iterations or is reciprocal of the desired accuracy. A  growing mini-batch, 
however, is undesirable in practice as requiring a large collection of samples per iteration  may easily prolong the duration of each 
iterate  without updating optimization parameters frequently enough \cite{defazio2018ineffectiveness}.   A 
notable exception to this trend is the the work of 
\cite{mokhtari2018stochastic} which employs a momentum variance-reduction 
technique requiring only one sample per iteration; however, this method suffers 
from suboptimal convergence rates. At the heart of this paper is the answer to the following question:
\vspace{-2mm}
\begin{quote}
	\textit{Can we achieve the optimal complexity bounds for a stochastic variant 
	of Frank-Wolfe while using a single stochastic sample per iteration?}
\end{quote}
\vspace{-2mm}
We show that the answer to the above question is positive and present the first projection-free method 
that  requires only one sample per iteration to update the optimization variable 
and yet  achieves the  optimal complexity bounds for convex, nonconvex, and monotone
DR-submodular settings.

More formally, we focus on a general \emph{non-oblivious} constrained 
stochastic optimization problem
\begin{equation}\label{eq:formation}
\min_{\bx \in \constraint}F(\bx) \triangleq \min_{x \in \constraint} 
\expect_{\bz \sim p(\bz;\bx)}[\tF(\bx;\bz)],
\end{equation}
where $\bx \in \mathbb{R}^d$ is the optimization variable, $\constraint 
\subseteq \mathbb{R}^d$ is the convex constraint set, and the objective 
function $F: \mathbb{R}^d \to \mathbb{R}$ is defined as the expectation over a 
set of functions $\tF$. The function $\tF: \mathbb{R}^d \times 
\mathcal{Z} \to \mathbb{R}$ is determined by $\bx$ and a 
random variable $\bz \in \mathcal{Z}$ with distribution $\bz \sim p(\bz;\bx)$. 
We refer to problem~\eqref{eq:formation} as a non-oblivious stochastic optimization problem 
as the 
distribution of the random variable $\bz$ depends on the choice of $\bx$. When 
the distribution $p$ is independent of $\bx$,  we are in the standard 
oblivious 
stochastic optimization regime where the goal is to solve
\begin{equation}\label{eq:ob_op}
\min_{\bx \in \constraint} F(\bx) \triangleq \min_{\bx \in \constraint} 
\expect_{\bz \sim p(\bz)}[\tF(\bx;\bz)].
\end{equation}
Hence, the oblivious problem  \eqref{eq:ob_op} can be considered as a special 
case of the non-oblivious problem  \eqref{eq:formation}. Note that 
non-oblivious stochastic optimization has broad applications in machine 
learning, including multi-linear extension of a discrete submodular 
function \citep{hassani2019stochastic}, MAP inference in determinantal
point processes (DPPs) \citep{kulesza2012determinantal}, and reinforcement 
learning 
\citep{sutton2018reinforcement,shen2019hessian}. 

%

Our goal is to propose an efficient FW-type method for the non-oblivious 
optimization problem  \eqref{eq:formation}. Here, the efficiency is measured 
by the number of 
stochastic oracle queries, {i.e.}, the sample complexity of $\bz$. 
As we mentioned earlier, among the stochastic variants of FW, the momentum 
stochastic Frank-Wolfe method proposed in 
\citep{mokhtari2017conditional,mokhtari2018stochastic} is the only method that 
requires only one sample per iteration. However, the stochastic oracle 
complexity of this  algorithm is suboptimal, {i.e.}, 
$\mathcal{O}(1/\epsilon^3)$ stochastic queries are required for both convex 
minimization and monotone DR-submodular maximization problems. This suboptimal 
rate is due to the fact that the gradient estimator in momentum FW is biased 
and it is necessary to use a more conservative averaging parameter to control the effect of the bias term.

\begin{table}[t!]
	\begin{center}
		\caption{Convergence guarantees of stochastic Frank-Wolfe methods for 
			constrained \textbf{convex} minimization}\label{table_convex}
		\begin{tabular}{|c|c|c|c|c|c|c|}
			\hline
			\textbf{Ref.}&   \textbf{batch}  &\!\textbf{complexity} & 
			\textbf{oblivious}& \textbf{non-oblivious}\\
			\hline
			\citep{DBLP:conf/icml/HazanK12} &  $\mathcal{O} (1/\epsilon^2)$  
			&$\mathcal{O}(1/ \epsilon^{4}) $ & \checkmark& \ding{55}\\
			\hline
			\citep{hazan2016variance} & $\mathcal{O} (1/\epsilon^2)$  
			&$\mathcal{O}(1/ \epsilon^{3}) $ & \checkmark& \ding{55}\\
			\hline
			\citep{mokhtari2018stochastic} &  $1$ &$\mathcal{O}(1/ 
			\epsilon^{3}) $ & \checkmark& \ding{55}\\
			\hline
			\citep{yurtsever2019conditional}&  $\mathcal{O}(1/\epsilon)$ 
			&$\mathcal{O}(1/\epsilon^2) $ & 
			\checkmark& \ding{55}\\
			\hline
			\citep{hassani2019stochastic}&  {{$\mathcal{O}(1/\epsilon)$}}  & 
			$\mathcal{O}(1/ \epsilon^{2}) $ & \checkmark& \checkmark\\
			\hline
			This paper &  1 & $\mathcal{O}(1/ \epsilon^{2}) $ & \checkmark& 
			\checkmark\\
			\hline
		\end{tabular}
		\vspace{-1mm}
		\vspace{-1mm}
	\end{center}
\end{table}

\begin{table}[t!]
	\begin{center}
		\caption{Convergence guarantees of stochastic Frank-Wolfe methods for 
			\textbf{non-convex} minimization}\label{table_nonconvex}
		\begin{tabular}{|c|c|c|c|c| c|c|}
			\hline
			\textbf{Ref.}&  \textbf{batch} & \textbf{complexity}  & 
			\textbf{oblivious}& \textbf{non-oblivious}\\
			\hline
			\citep{hazan2016variance} &  $\mathcal{O} (1/\epsilon^2)$ 
			&$\mathcal{O}(1/ \epsilon^{4}) $ & \checkmark& \ding{55}\\
			\hline
			\citep{hazan2016variance}  & $\!\!\!\mathcal{O} 
			(1/\epsilon^{4/3})\!\!\!$ &$\mathcal{O}(1/ \epsilon^{10/3}) $ & 
			\checkmark& \ding{55}\\
			\hline
			\citep{shen2019complexities} & {{$\mathcal{O}(1/\epsilon)$}}   
			&$\mathcal{O}(1/ \epsilon^{3}) $ & \checkmark& \ding{55}\\
			\hline
			\citep{yurtsever2019conditional}&  {{$\mathcal{O}(1/\epsilon)$}}  
			&$\mathcal{O}(1/\epsilon^3) $ & 
			\checkmark& \ding{55}\\
			\hline
			\citep{hassani2019stochastic} &  {{$\mathcal{O}(1/\epsilon)$}}   
			&$\mathcal{O}(1/ \epsilon^{3}) $ & \checkmark& \checkmark	\\
			\hline
			This paper & 1 & $\mathcal{O}(1/ \epsilon^{3}) $ & \checkmark& 
			\checkmark\\
			\hline
		\end{tabular}
		\vspace{-1mm}
		\vspace{-1mm}
	\end{center}
\end{table}

\begin{table}[t!]
	\begin{center}
		\caption{Convergence guarantees for stochastic \textbf{monotone 
				continuous DR-submodular} function 
			maximization}\label{table_sub}
		\begin{tabular}{| c| c| c| c| c| c| }
			\hline
			\textbf{Ref.}&   \textbf{batch} &\textbf{utility} & 
			\textbf{complexity}  \\
			\hline 
			\citep{hassani2017gradient} & 1 &$(1/2)\rm{OPT}$$-\epsilon$& 
			$O(1/{\epsilon^{2}})$ \\
			\hline
			\citep{mokhtari2018stochastic}  & 1 &$(1-1/e)\rm{OPT}$$-\epsilon$& 
			$O(1/{\epsilon^{3}})$ \\
			\hline
			\citep{hassani2019stochastic}  &{{$\mathcal{O}(1/\epsilon)$}} 
			&$(1-1/e)\rm{OPT}$$-\epsilon$& 
			{\color{black}{$\mathcal{O}(1/\epsilon^2)$}} \\
			\hline
			This paper & 1 &  $(1-1/e)\rm{OPT}$$-\epsilon$& $\mathcal{O}(1/ 
			\epsilon^{2}) $ \\
			\hline
		\end{tabular}
		\vspace{-1mm}
		\vspace{-1mm}
	\end{center}
\end{table}

Theoretical results of \sfw and other related works are 
summarized in  Tables~\ref{table_convex}-\ref{table_sub}. These results show 
that  \sfw attains the optimal or best known complexity bounds in all the 
considered settings, while requiring only \emph{one single} stochastic oracle 
query per iteration and avoiding large batch sizes altogether. 

To resolve this issue, we propose a one-sample stochastic Frank-Wolfe method, 
called \sfw, which modifies the gradient approximation in momentum FW to ensure 
that the resulting gradient estimation is an unbiased estimator of the gradient 
(Section~\ref{sec:sfw}). This goal has been achieved by adding an unbiased 
estimator of the gradient variation  $\Delta_t = \nabla F(\bx_t) - \nabla 
F(\bx_{t-1})$ to the gradient approximation vector 
(Section~\ref{sec:sg_approx}). We later explain why coming up with an unbiased 
estimator of the gradient difference $\Delta_t$ could be a challenging task in the non-oblivious setting and 
show how we overcome this difficulty (Section~\ref{sec:grad_var_est}). We also 
characterize the convergence guarantees of \sfw for convex minimization, 
nonconvex minimization and monotone DR-submodular maximization 
(Section~\ref{sec:results}). In particular, we show that \sfw achieves the 
optimal  convergence rate  of $\mathcal{O}(1/\epsilon^2)$  for reaching an  
$\epsilon$-suboptimal solution in the stochastic convex setting, and  a 
$(1-1/e)-\epsilon$ approximate solution for a stochastic monotone DR-submodular 
maximization problem. Moreover, in a general non-convex setting, \sfw finds an 
$\epsilon$-first-order stationary point after at most 
$\mathcal{O}(1/\epsilon^3)$ iterations, achieving the current best known 
convergence rate.  Finally, we study the oblivious problem in \eqref{eq:ob_op} 
and show that our proposed \sfw method becomes significantly simpler and  the
corresponding theoretical results hold under  less strict assumptions. We further highlight 
the similarities between the variance reduced method in 
\citep{cutkosky2019momentum} also known as STORM and the oblivious variant of 
\sfw. Indeed, our algorithm has been originally inspired by STORM.

\section{Related Work}
As a projection-free algorithm, Frank-Wolfe method \citep{frank1956algorithm} 
has been studied for both convex optimization 
\citep{jaggi2013revisiting,lacoste2015global,garber2015faster,hazan2016variance,mokhtari2018stochastic}
and non-convex optimization problems
\citep{lacoste2016convergence,reddi2016stochastic,mokhtari2018escaping,shen2019hessian,hassani2019stochastic}.
In large-scale settings, distributed FW methods were proposed to solve 
specific 
problems, including optimization under block-separable constraint set 
\citep{wang2016parallel}, and learning low-rank matrices 
\citep{zheng2018distributed}. The communication-efficient distributed FW 
variants were proposed for specific sparse learning problems in 
\cite{bellet2015distributed,lafond2016d}, and for general constrained 
optimization problems in \citep{zhang2019quantized}. Zeroth-order FW methods 
were studied in \citep{sahu2018towards,chen2019black}.

Several works have studied different ideas for reducing variance in stochastic 
settings.  The SVRG method was proposed by
\cite{Johnson2013Accelerating} for the convex setting and then extended to the 
nonconvex setting by several other works 
\citep{allen2016variance,reddi2016stochastic,zhou2018stochastic}. The StochAstic
Recursive grAdient algoritHm (SARAH)  was studied in 
\citep{nguyen2017sarah,nguyen2017stochastic}. Then as a variant of SARAH, the 
Stochastic Path-Integrated Differential Estimator (SPIDER) technique was 
proposed by \cite{fang2018spider}. Based on SPIDER, various algorithm for 
convex and non-convex optimization problems have been studied 
\citep{shen2019complexities,hassani2019stochastic,yurtsever2019conditional}. 

In this paper, we also consider optimizing  an important  subclass of non-convex objectives, known as continuous DR-submodular 
functions that generalize the diminishing returns property to the continuous domains. 
%
 Continuous 
DR-submodular functions can be minimized exactly 
\citep{bach2015submodular,staib2017robust}, and maximized approximately 
\citep{bian16guaranteed,bian2017continuous,hassani2017gradient,mokhtari2017conditional}.
They have interesting applications in machine learning, including 
experimental design 
\citep{Chen2018Online}, MAP inference in determinantal
point processes (DPPs) \citep{kulesza2012determinantal}, and mean-field 
inference in probabilistic models \citep{bian2018optimal}.

\section{One Sample SFW Algorithm}\label{sec:sfw}

In this section, we introduce our proposed one sample SFW (\sfw) method. We 
first present the mechanism for computing a variance reduced unbiased estimator 
of the gradient $\nabla F(\bx_t)$. Then, we explain the procedure for computing 
an unbiased estimator of the gradient variation $\Delta_t = \nabla F(\bx_t) - 
\nabla F(\bx_{t-1})$ in a non-oblivious setting which is required for the 
gradient approximation of \sfw. Then, we present the complete description of 
our proposed method.

\subsection{Stochastic gradient approximation}\label{sec:sg_approx}


In our work, we build on the momentum variance 
reduction approach proposed in 
\citep{mokhtari2017conditional,mokhtari2018stochastic}
to reduce the 
variance of the one-sample method. 
To be more precise, in the momentum FW method 
\citep{mokhtari2017conditional}, we update the gradient approximation $\bd_t$ 
at round 
$t$ according to the update
\begin{equation}\label{eq:momentum}
\bd_t = (1-\rho_t) \bd_{t-1} + \rho_t \nabla \tF(\bx_t; \bz_t),
\end{equation}
where $\rho_t$ is the averaging parameter and $\nabla \tF(\bx_t; \bz_t)$ is a 
\emph{one-sample} estimation of the gradient. Since $\bd_t$ is a weighted 
average of the previous gradient estimation $\bd_{t-1}$ and the newly updated 
stochastic gradient, it has a lower variance comparing to one-sample 
estimation $\nabla \tF(\bx_t; \bz_t)$. In particular, it was shown by 
\citet{mokhtari2017conditional} that the variance of gradient approximation in 
\eqref{eq:momentum} approaches zero at a sublinear rate of $O(t^{-2/3})$. The 
momentum approach reduces the variance of gradient approximation, but it leads 
to a \emph{biased} gradient approximation, i.e., $\bd_t$  is not an 
unbiased estimator of the gradient $\nabla F(\bx_t)$. Consequently, it is 
necessary to use a conservative averaging parameter $\rho_t$ for momentum FW to 
control the effect of the bias term which leads to a 
sublinear error rate of $\mathcal{O}(t^{-1/3})$ and overall complexity of 
$\mathcal{O}(1/\epsilon^3)$.


To resolve this issue and come up with a fast momentum based FW method for the 
non-oblivious problem in \eqref{eq:formation}, we slightly modify the gradient 
estimation in \eqref{eq:momentum} to ensure that the resulting gradient 
estimation is an unbiased estimator of the gradient $\nabla F(\bx_t)$. 
Specifically, we add the term  
$\tilde{\Delta}_t$, which is an unbiased estimator of the gradient variation  
$\Delta_t = \nabla F(\bx_t) - \nabla F(\bx_{t-1})$, to $\bd_{t-1}$. This 
modification leads to the following gradient approximation
\begin{equation}\label{eq:momentum_new}
\bd_t = (1-\rho_t) (\bd_{t-1}+\tilde{\Delta}_t) + \rho_t \nabla \tF(\bx_t; 
\bz_t).
\end{equation}
To verify that $\bd_t $ is an unbiased estimator of $\nabla F(\bx_t)$ we can  
use a simple induction argument. Assuming that $\bd_{t-1}$ is an unbiased 
estimator of $\nabla F(\bx_t)$ and $\tilde{\Delta}_t$ is an unbiased estimator 
of $\nabla F(\bx_t) - \nabla F(\bx_{t-1})$ we have $\expect [\bd_t] = 
(1-\rho_t)( \nabla F(\bx_{t-1})+(\nabla F(\bx_t) - \nabla F(\bx_{t-1})))+\rho_t 
\nabla F(\bx_t)= \nabla F(\bx_t)$. Hence, the gradient approximation in 
\eqref{eq:momentum_new} leads to an unbiased approximation of the gradient. Let 
us now explain how to compute an unbiased estimator of the gradient variation 
$\Delta_t = \nabla F(\bx_t) - \nabla F(\bx_{t-1})$ for a non-oblivious setting. 

\subsection{Gradient variation estimation}\label{sec:grad_var_est}

The most natural approach for estimating the gradient variation  
$\Delta_t=\nabla F(\bx_t) - \nabla F(\bx_{t-1})$ using only one sample $\bz$ is 
computing the difference of two consecutive stochastic gradients, i.e., 
$\nabla \tF(\bx_t; \bz)-\nabla \tF(\bx_{t-1}; \bz)$. However, this approach 
leads to an unbiased estimator of the gradient variation  $\Delta_t$ only in 
the oblivious setting where $p(\bz)$ is independent of the choice of $\bx$, and 
would introduce bias in the more general non-oblivious case.  To better 
highlight this issue, assume that $\bz$ is sampled according to distribution 
$p(\bz;\bx_t)$. Note that $\nabla \tF(\bx_t; \bz)$ is an unbiased estimator of 
$\nabla F(\bx_t)$, i.e., $\expect[\nabla F(\bx_t; \bz)] = \nabla 
F(\bx_t)$, 
however, $\nabla \tF(\bx_{t-1};\bz)$ is not an unbiased estimator of $\nabla 
F(\bx_{t-1})$ since $p(\bz;\bx_{t-1})$ may be different from $p(\bz;\bx_t)$.

To circumvent this obstacle, an \emph{unbiased} estimator of $\Delta_t$ was 
introduced in \citet{hassani2019stochastic}. To explain their proposal for 
approximating the 
gradient variation using only one sample, note that the difference 
$\Delta_t=\nabla F(\bx_t) - \nabla F(\bx_{t-1})$ can be written as
\begin{equation*}
\Delta_t = \int_0^1 \nabla ^2 F(\bx_t(a))(\bx_t -\bx_{t-1}) \mathrm{d}a = 
\left[\int_0^1 \nabla ^2 F(\bx_t(a))\mathrm{d}a \right] (\bx_t -\bx_{t-1}),
\end{equation*}
where $\bx_t(a) = a\bx_t + (1-a)\bx_{t-1}$ for  $a\in [0,1]$. According to this 
expression, one can find an unbiased estimator of $\int_0^1 \nabla 
^2 F(\bx_t(a))\mathrm{d}a$ and use its product with $(\bx_t -\bx_{t-1})$ to 
find an unbiased estimator of $\Delta_t$. It can be easily verified that 
$\nabla ^2 F(\bx_t(a))(\bx_t -\bx_{t-1})$ is an unbiased estimator of 
$\Delta_t$ if $a$ is chosen from $[0,1]$ uniformly at random. Therefore, all we 
need is to come up with an unbiased estimator of the Hessian $\nabla ^2 F$.

By basic calculus, we can show that  $\forall \bx \in 
\constraint$ and $\bz$ with distribution $p(\bz;\bx)$, the matrix 
$\tilde{\nabla}^2F(\bx;\bz)$ defined as
\begin{align}
\tilde{\nabla}^2F(\bx;\bz) &= \tF(\bx;\bz)[\nabla \log p(\bz;\bx)][\nabla \log 
p(\bz;\bx)]^\top + \nabla^2 \tF(\bx;\bz)  + [\nabla \tF (\bx; 
\bz)][\nabla \log 
p(\bz;\bx)]^\top \nonumber \\
&\ + \tF(\bx;\bz) \nabla^2 \log p(\bz;\bx)  + [\nabla \log 
p(\bz;\bx)][\nabla \tF(\bx;\bz)]^\top,
\label{eq:gradient_diff_estimation2}
\end{align}
is an \emph{unbiased} estimator of $\nabla^2F(\bx)$.  Note that the above 
expression requires only one sample of $\bz$. As a result, we can construct 
$\tilde{\Delta}_t$ as an unbiased estimator of ${\Delta}_t$ using only one 
sample
\begin{equation}\label{eq:gradient_diff_estimation}
\tilde{\Delta}_t \triangleq \tilde{\nabla}_t^2 (\bx_t - \bx_{t-1}),
\end{equation}
where $\tilde{\nabla}_t^2 = \tilde{\nabla}^2F(\bx_t(a);\bz_t(a))$,
 and $\bz_t(a)$ follows the distribution $p(\bz_t(a);\bx_t(a))$. By using this 
procedure, we can indeed compute the vector $\bd_t$ in \eqref{eq:momentum_new} 
with only one sample of $\bz$ per iteration.

We note that if we only use one sample of $\bz$ per iteration, i.e., 
$\bz_t = \bz_t(a)$, the 
stochastic gradient in \eqref{eq:momentum} becomes $\nabla \tF(\bx_t; 
\bz_t(a))$. It is not an unbiased estimator of $\nabla F(\bx_t)$, since 
$\nabla F(\bx_t) =\expect_{\bz\sim p(\bz;\bx_t)}\tF(\bx_t;\bz)$, while 
$\bz_t(a) \sim p(\bz; \bx_t(a))$. Through a completely
different analysis from the ones in 
\citep{mokhtari2017conditional,hassani2019stochastic}, we 
show that the modified $\bd_t$ is still a good gradient estimation 
(\cref{lem:graident_error}), which allows
the establishment of the optimal stochastic oracle
complexity for our proposed algorithm.

Another issue of this scheme is that in 
\eqref{eq:gradient_diff_estimation2} and \eqref{eq:gradient_diff_estimation}, 
we need to 
calculate $\nabla^2 \tF(\bx_t(a);\bz_t(a))(\bx_t-\bx_{t-1})$ and $\nabla^2 \log 
p(\bx_t(a);\bz_t(a))(\bx_t-\bx_{t-1})$, where computation of Hessian is 
involved. When exact Hessian is not accessible, however, we can resort to an 
approximation by the difference of two gradients.

Precisely, for any function $\psi: \mathbb{R}^d \to \mathbb{R}$, any vector 
$\bu \in \mathbb{R}^d$ with $\|\bu \|\leq D$, and some $\delta >0$ small 
enough, we have
\begin{equation*}
\phi(\delta;\psi) \triangleq \frac{\nabla \psi(\bx+\delta \bu) - \nabla \psi(x 
- \delta\bu)}{2\delta} \approx \nabla^2 \psi(\bx) \bu.
\end{equation*}
If we assume that $\psi$ is $L_2$-second-order smooth, i.e., 
$\|\nabla^2\psi(\bx) - \nabla^2\psi(\by)) \| \leq L_2 \|\bx-\by \|,\ \forall 
\bx, \by \in \mathbb{R}^d$, we can 
upper bound the 
approximation error quantitatively:
\begin{align}\label{eq:gradient_difference_bound}
\|\nabla^2 \psi(\bx) \bu - \phi(\delta;\psi) \| = \|\nabla^2 \psi(\bx) \bu - 
\nabla^2 \psi(\tilde{\bx}) \bu) \| \leq D^2L_2\delta,
\end{align}
where $\tilde{\bx}$ is obtained by the mean-value theorem. In other words, the 
approximation error can be sufficiently small for proper $\delta$. So we can 
estimate $\Delta_t$ by
\begin{align}\label{eq:gradient_diff_estimation3}
\tilde{\Delta}_t &= \tF(\bx;\bz)[\nabla \log p(\bz;\bx)][\nabla \log 
p(\bz;\bx)]^\top \bu_t + \phi(\delta_t,\tF(\bx;\bz))  + [\nabla \tF (\bx; 
\bz)][\nabla \log 
p(\bz;\bx)]^\top\bu_t \nonumber\\
&\ + \tF(\bx;\bz) \phi(\delta_t,\log p(\bz,\bx))  + [\nabla \log 
p(\bz;\bx)][\nabla \tF(\bx;\bz)]^\top \bu_t,
\end{align}
where $\bu_t = \bx_t -\bx_{t-1}$, $\bx, \bz, \delta_t$ are chosen 
appropriately. We also 
note that since computation of gradient difference has a 
complexity of $\mathcal{O}(d)$, while that for Hessian is $\mathcal{O}(d^2)$, 
this approximation strategy can also help to accelerate the optimization 
process.

\subsection{Variable update}

Once the gradient approximation $\bd_t$ is computed, we can follow the update 
of conditional gradient methods for computing the iterate $\bx_t$. In this 
section, we introduce two different schemes for updating the iterates depending 
on the problem that we aim to solve. 

For minimizing a general (non-)convex function using one sample stochastic FW, 
we update the iterates according to the update 
\begin{equation}\label{eq:option1}
\bx_{t+1} = \bx_{t}+\eta_{t} (\bv_t-\bx_{t}),
\text{ where } \bv_t = 
\argmin_{v \in \constraint}\{\bv^\top \bd_{t}\}.
\end{equation} 
In this case, we find the direction that minimizes the inner product with the 
current gradient approximation $\bd_t$ over the constraint set $\constraint$, 
and the updated variable $\bx_{t+1}$ by descending in the direction of 
$\bv_t-\bx_{t}$ with step size $\eta_t$.

For monotone DR-submodular maximization, the update rule is slightly different, 
and a stochastic variant of the continuous greedy method 
\citep{vondrak2008optimal} can be used. 
Using the same stochastic estimator $\bd_t$ as in the (non-)convex case, the 
update rule for DR-Submodular optimization is given by
\begin{equation}\label{eq:option2}
\bx_{t+1} = \bx_{t}+\eta_{t} \bv_t, \quad
\text{where } \bv_t = 
\argmax_{v \in \constraint}\{\bv^\top \bd_{t}\}
\end{equation}
where $\eta_t=1/T$. Note that in this case we find the direction that maximizes 
the inner product with the current gradient approximation $\bd_t$ over the 
constraint set $\constraint$, and move towards that direction with step size 
$\eta_t=1/T$. Hence, if we start from the origin, after $T$ steps the outcome 
will be a feasible point as it can be written as the average of $T$ feasible 
points. 

The description of  our proposed \sfw method for smooth (non-)convex 
minimization as well as monotone DR-submodular maximization is outlined in 
\eqref{alg:one_sample}.

\begin{algorithm}[t]
	\caption{One-Sample SFW}
	\begin{algorithmic}[1]
		\Require Step sizes $\rho_t \in 
		(0,1),   
		\eta_t \in (0,1)$, initial point $\bx_1 \in 
		\constraint$, total number 
		of iterations $T$
		\Ensure $\bx_{T+1}$ or $\bx_o$, where $\bx_o$ is chosen from $\{\bx_1, 
		\bx_2,\cdots, \bx_T\}$ uniformly at random  
		\For{$t=1,2,\dots, T$}
		\If{$t=1$} 
		\State Sample a point $\bz_1$ according 
		to $p(\bz_1,\bx_1)$
		\State  Compute $\bd_1 = \nabla \tilde{F}(\bx_1; 
		\bz_1)$		
		\Else
		\State Choose $a$ uniformly at random from $[0,1]$
		\State Compute $\bx_t(a) = a\bx_t + (1-a)\bx_{t-1}$
		\State Sample a point $\bz_t$ according to $p(\bz; \bx_t(a))$
		\State Compute $\tilde{\Delta}_t$ either by $\tnabla_t^2 = \tnabla^2 
		F(\bx_t(a); \bz_t)$ based on \eqref{eq:gradient_diff_estimation2} and 
		$\tilde{\Delta}_t =\tilde{\nabla}_t^2 (\bx_t - \bx_{t-1})$ (Exact 
		Hessian Option); or by \cref{eq:gradient_diff_estimation3} with 
		$\bx=\bx_t(a), \bz=\bz_t$ (Gradient Difference Option)
		\State $\bd_t = (1-\rho_t)(\bd_{t-1}+\tilde{\Delta}_t) + \rho_t \nabla 
		\tilde{F}(\bx_t,\bz_t)$
		\EndIf
		\State (non-)convex min.: Update $\bx_{t+1}$ based on 
		\eqref{eq:option1} 
		\State  DR-sub. max.: Update $\bx_{t+1}$ based on \eqref{eq:option2} 
		\EndFor
	\end{algorithmic}
	\label{alg:one_sample}
\end{algorithm}

\section{Main Results}\label{sec:results}
Before presenting the convergence results of our algorithm, we first state our 
assumptions on the constraint set $\constraint$, the 
stochastic function $\tF$, and the distribution $p(\bz;\bx)$. 

\begin{assump}\label{assum:constraint}
	The constraint set $\constraint \subseteq \mathbb{R}^d$ is compact with 
	diameter $D = max_{x,y \in \constraint}\|x-y\|$, and radius $R = 
	\max_{x\in\constraint} \|x\|$. 
\end{assump}


\begin{assump}\label{assum:stoch_bound}
The stochastic function $\tF(\bx;\bz)$ has uniformly bounded
function value, i.e., $| \tF(\bx;\bz)| \leq B$ for all $\bx \in 
\constraint,  \bz \in \mathcal{Z}$. 
\end{assump}

\begin{assump}\label{assum:gradient_norm}
The stochastic gradient $\nabla \tF$ has uniformly bound norm, \emph{i.e.}, 
$\|\nabla 
\tF(\bx;\bz) \| \leq G_{\tF}, \forall \bx \in \constraint, \forall \bz \in 
\mathcal{Z}$. The norm of the gradient of $\log p$ has bounded fourth-order 
moment, \emph{i.e.}, $\expect_{\bz \sim p(\bz;\bx)} \|\nabla \log p(\bz;\bx) 
\|^4 \leq 
G_p^4$. We also define $G=\max\{G_\tF, G_p\}$.
\end{assump}

\begin{assump}\label{assum:second_order}
The stochastic Hessian $\nabla^2 \tF$ has uniformly bounded spectral norm: $\| 
\nabla^2 \tF(\bx;\bz)\| \leq L_{\tF}, \forall \bx \in \constraint, \forall 
\bz \in \mathcal{Z}$. The spectral norm of the Hessian of $\log p$ has bounded 
second-order moment: $\expect_{\bz \sim p(\bz;\bx)} \|\nabla^2 \log p(\bz;\bx) 
\|^2 \leq L_p^2$. We also define $L = \max \{L_\tF, L_p\}$. 
\end{assump}

We note that in Assumptions~\ref{assum:stoch_bound}-\ref{assum:second_order}, 
we 
assume that the stochastic function $\tF$ has uniformly bounded function 
value, gradient norm, and second-order differential. Moreover, with these 
assumptions,  we can establish an upper bound for the second-order moment of 
the spectral norm of the Hessian estimator $\tnabla^2 F(\bx;\bz)$, which is 
defined 
in \eqref{eq:gradient_diff_estimation2}.

\begin{lemma}\label{lem:scg++8.1}[Lemma 7.1 of 
	\citep{hassani2019stochastic}]
Under Assumptions~\ref{assum:stoch_bound}-\ref{assum:second_order}, for all 
$\bx \in \constraint$, we have
\begin{equation*}
\expect_{\bz \sim p(\bz;\bx)}[\| \tnabla^2F(\bx;\bz) \|^2] \\
\leq{} 4B^2G^4 + 
16G^4 + 4L^2 + 4B^2L^2 
\triangleq{} \bar{L}.
\end{equation*}
\end{lemma}
Note that the result in Lemma \eqref{lem:scg++8.1} also implies the 
$\bar{L}$-smoothness of $F$, since 
\begin{align*}
\|\nabla^2 
F(\bx) \|^2 = \|\expect_{\bz \sim p(\bz;\bx)}[\tnabla^2 F(\bx;\bz)] \|^2
\leq \expect_{\bz \sim p(\bz;\bx)}[\|\tnabla^2 F(\bx;\bz) \|^2] \leq \bar{L}^2.
\end{align*}
In other words, the conditions in 
Assumptions~\ref{assum:stoch_bound}-\ref{assum:second_order} implicitly imply 
that the objective function $F$ is $\bar{L}$-smooth.

To establish the convergence guarantees for our proposed \sfw algorithm, the 
key step is to derive an upper bound on the errors of the estimated gradients. 
To do so, we prove the following lemma, which provides the required upper 
bounds in 
different settings of parameters.

\begin{lemma}\label{lem:graident_error}
Consider the gradient approximation $\bd_t$ defined in \eqref{eq:momentum_new}. 
Under Assumptions~\ref{assum:constraint}-\ref{assum:second_order},
 if we run Algorithm~\ref{alg:one_sample} with Exact Hessian Option in Line 9, 
 and with parameters $\rho_t = (t-1)^{-\alpha}, \forall t \ge 2$, and 
$\eta_t \leq t^{-\alpha}, \forall t \ge 1$ for some $\alpha \in (0,1]$, 
then the gradient estimation $\bd_t$ satisfies
\begin{equation}
\expect[\|\nabla F(\bx_t) - \bd_t\|^2] \leq C t^{-\alpha},
\end{equation}
where the constant $C$ is given by
 $$C\! =\! \max 
 \!\left\{\!\frac{2(2G\!+\!D\bar{L})^2}{2\!-\!2^{-\alpha}\!-\!\alpha}, 
\left[\frac{2}{2\!-\!2^{-\alpha}\!-\!\alpha} \right]^4\!\!, 
[2D(\bar{L}\!+\!L)]^4 \!\right\}\!.$$	
\end{lemma}

Lemma \eqref{lem:graident_error} shows that with appropriate parameter setting, 
the 
gradient error converges to zero at a rate of $t^{-\alpha}$. With this 
unifying upper bound, we can obtain the convergence rates of our algorithm for 
different kinds of objective functions. 

If in the update of \sfw we use the Gradient Difference Option in Line 9 of 
\cref{alg:one_sample} to estimate $\tilde{\Delta}_t$, as pointed out above, we 
need one further assumption on second-order smoothness of the functions $\tF$ 
and $\log p$. 
\begin{assump}\label{assum:hessian}
The stochastic function $\tF$ is uniformly $L_{2,\tF}$- second-order smooth: 
$\|\nabla^2 \tF(\bx;\bz) - 
\nabla^2 \tF(\by;\bz) \|\leq L_{2,\tF}\|\bx-\by \|,\ \forall \bx,\by \in 
\constraint, \forall \bz \in \mathcal{Z}$. The log probability $\log 
p(\bz;\bx)$ is uniformly $L_{2,p}$-second-order smooth: $\|\nabla^2 \log 
p(\bz;\bx) - 
\nabla^2 \log p(\bz;\by) \|\leq L_{2,p}\|\bx-\by \|,\ \forall \bx,\by \in 
\constraint, \forall \bz \in \mathcal{Z}$. We also define $L_2 = 
\max\{L_{2,\tF}, L_{2,p}\}$.
\end{assump}

We note that under \eqref{assum:hessian}, the approximation bound in 
\eqref{eq:gradient_difference_bound} holds for both 
$\tF$ and $\log p$. So for $\delta_t$ sufficiently small, the error introduced 
by the Hessian approximation can be ignored. Thus similar upper bound for 
errors of estimated gradient still holds. 
\begin{lemma}\label{lem:graident_error_2}
	Consider the gradient approximation $\bd_t$ defined in 
	\eqref{eq:momentum_new}. Under 
	Assumptions~\ref{assum:constraint}-\ref{assum:hessian},
	if we run Algorithm~\ref{alg:one_sample} with Gradient Difference Option in 
	Line 9, and with parameters $\rho_t = (t-1)^{-\alpha}, \delta_t = 
	\frac{\sqrt{3}\eta_{t-1\bar{L}}}{DL_2(1+B)},\ \forall t \ge 2$, and 
	$\eta_t \leq t^{-\alpha}, \forall t \ge 1$ for some $\alpha \in (0,1]$, 
	then the gradient estimation $\bd_t$ satisfies
	\begin{equation}
	\expect[\|\nabla F(\bx_t) - \bd_t\|^2] \leq C t^{-\alpha},
	\end{equation}
	where the constant $C$ is given by
	\begin{equation*}
	C= \max \
	\!\bigg\{\!\frac{8(D^2\bar{L}^2+G^2+GD\bar{L})}{2\!-\!2^{-\alpha}\!-\!\alpha},
	\left(\frac{2}{2\!-\!2^{-\alpha}\!-\!\alpha} \right)^4\!\!, \quad 
	[4D(\bar{L}\!+\!L)]^4 \bigg\}.
	\end{equation*}	
\end{lemma}
\cref{lem:graident_error_2} shows that with Gradient Difference Optionin Line 9 
of \cref{alg:one_sample}, the error of estimated gradient can obtain the same 
order of convergence rate as that with Exact Hessian Option. So in the 
following three subsections, we will present the theoretical results of our 
proposed \sfw 
algorithm with Exact Hessian Option, for convex minimization, non-convex 
minimization, and monoton DR-submodular maximization, respectively. The results 
of Gradient Difference Option only differ in a factor of constant.

\subsection{Convex Minimization}
For convex minimization problems, to obtain an $\epsilon$-suboptimal solution, 
\eqref{alg:one_sample} only requires at most $\mathcal{O}(1/\epsilon^2)$ 
stochastic oracle queries, and $\mathcal{O}(1/\epsilon^2)$ linear optimization 
oracle calls. Or precisely, we have 
\begin{theorem}[Convex]\label{thm:convex}
Consider the \sfw method outlined in Algorithm~\ref{alg:one_sample}  
with Exact Hessian Option in Line 9. Further, suppose the conditions in 
Assumptions~\ref{assum:constraint}-\ref{assum:second_order} hold, and  assume 
that $F$ is convex on $\constraint$. Then, if we set the algorithm parameters 
as $\rho_t = 
 (t-1)^{-1}$ and $\eta_t = t^{-1}$, then the output is feasible $\bx_{T+1} \in 
 \constraint$ and satisfies
\begin{equation*}
\expect[F(\bx_{T+1}) - F(\bx^*)] \leq \frac{2\sqrt{C}D}{\sqrt{T}} + 
\frac{\bar{L}D^2(1+\ln T)}{2T} ,
\end{equation*}
where $C=\max \{4(2G+D\bar{L})^2, 256, [2D(\bar{L}+L)]^4 \}$, and $\bx^*$ is a 
minimizer of $F$ on $\constraint$. 	
\end{theorem}

The result in Theorem~\ref{thm:convex} shows that the proposed one sample 
stochastic Frank-Wolfe method, in the convex setting, has an overall complexity 
of $\mathcal{O}(1/\epsilon^2)$ for finding an $\epsilon$-suboptimal solution. 
Note that to prove this claim we used the result in 
Lemma~\ref{lem:graident_error} for the case that $\alpha=1$, i.e., the 
variance of gradient approximation converges to zero at a rate of 
$\mathcal{O}(1/t)$.

\subsection{Non-Convex Minimization}
For non-convex minimization problems, showing that the gradient norm approaches 
zero, i.e., $\|\nabla F(\bx_t) \| \to 0$, implies convergence to a 
stationary point in the \emph{unconstrained} setting. Thus, it is usually used 
as a measure for convergence. In the constrained setting, however, the norm of 
gradient is not a proper measure for defining stationarity and we instead used 
the Frank-Wolfe Gap 
\citep{jaggi2013revisiting,lacoste2016convergence}, which is defined by
\begin{equation*}
\mathcal{G}(\bx) = \max_{\bv \in \constraint} \langle\bv-\bx,-\nabla F(\bx) 
\rangle.
\end{equation*}
We note that by definition, $\mathcal{G}(\bx) \geq 0,  \forall \bx \in 
\constraint$. If some point $\bx \in \constraint$ satisfies $\mathcal{G}(\bx) = 
0$, then it is a first-order stationary point.

In the following theorem, we formally prove the number of iterations required 
for one sample stochastic FW to find an $\epsilon$-first-order stationary point 
in expectation, i.e., a point $\bx$ that satisfies 
$\expect[\mathcal{G}(\bx)] \leq \epsilon$. 

\begin{theorem}[Non-Convex]\label{thm:nonconvex}
Consider the \sfw method outlined in Algorithm~\ref{alg:one_sample}  
with Exact Hessian Option in Line 9. Further, suppose the conditions in 
Assumptions~\ref{assum:constraint}-\ref{assum:second_order} hold. Then, if we 
set the algorithm parameters as 
$\rho_t 
= (t-1)^{-2/3},$ and $\eta_t = T^{-2/3}$, then the 
output is feasible $\bx_{o} \in \constraint$ and satisfies
\begin{equation*}
\expect[\mathcal{G}(\bx_o)] \leq \frac{2B+3\sqrt{C}D/2}{T^{1/3}} + 
\frac{\bar{L}D^2}{2T^{2/3}},
\end{equation*}
where the constant $C$ is given by
$$C\! = \!\max 
\left\{\frac{2(2G\!+\!D\bar{L})^2}{\frac{4}{3}-2^{-\frac{2}{3}}}, \left[ 
\frac{2}{\frac{4}{3}-2^{-\frac{2}{3}}}\right]^4\!\!, [2D(\bar{L}\!+\!L)]^4 
\right\}\!.$$
\end{theorem}
We remark that Theorem~\eqref{thm:nonconvex} shows that 
Algorithm~\ref{alg:one_sample} finds an 
$\epsilon$-first order stationary points after at most 
$\mathcal{O}(1/\epsilon^3)$ iterations, while uses exactly one stochastic 
gradient per iteration. Note that to obtain the best performance guarantee in 
Theorem~\eqref{thm:nonconvex}, we used the result of 
Lemma~\ref{lem:graident_error} for the case that $\alpha=2/3$, i.e., the 
variance of gradient approximation converges to zero at a rate of 
$\mathcal{O}(T^{-2/3})$.

\subsection{Monotone DR-Submodular Maximization}

In this section, we focus on the convergence properties of one-sample 
stochastic Frank-Wolfe or one-sample stochastic Continuous Greedy for solving a 
monotone DR-submodular maximization problem.
%
Consider a 
differentiable function $F: \mathcal{X} \to \mathbb{R}_{\geq 0}$, where the 
domain 
$\mathcal{X} \triangleq \prod_{i=1}^{d} \mathcal{X}_i$, and each 
$\mathcal{X}_i$ is a compact subset of $\mathbb{R}_{\geq 0}$. We say $F$ is 
continuous DR-submodular if for all $\bx, \by \in \mathcal{X}$ that satisfy 
$\bx \leq \by$ and every $i \in \{1,2,\cdots,d\}$, we have $\frac{\partial 
F}{\partial x_i}(\bx) \geq \frac{\partial F}{\partial x_i}(\by)$.

An important property of continuous DR-submodular function is the concavity 
along the 
non-negative directions \citep{calinescu2011maximizing,bian16guaranteed}: for 
all $\bx, \by \in \mathcal{X}$ such that $\bx \leq \by$, we have $F(\by) \leq 
F(\bx) + \langle \nabla F(\bx), 
\by-\bx \rangle$. We say $F$ is monotone if for all $\bx, \by \in 
\mathcal{X}$ such that $\bx \leq \by$, we have $F(\bx) \leq F(\by)$. 

For continuous DR-submodular maximization, it has been shown that approximated 
solution within a factor of $(1-e^{-1} + \epsilon)$ can not be obtained in 
polynomial time \citep{bian16guaranteed}. As a result, we analyze 
the convergence rate for the tight $(1-e^{-1})$OPT approximation. To achieve a 
$(1-e^{-1})\mathrm{OPT}- \epsilon$ approximation guarantee, our proposed 
algorithm requires at most $\mathcal{O}(1/\epsilon^2)$ 
stochastic oracle queries, and $\mathcal{O}(1/\epsilon^2)$ linear optimization 
oracle calls, as we show in the following theorem.
 
\begin{theorem}[Submodular]\label{thm:sub}
Consider the \sfw method outlined in Algorithm~\ref{alg:one_sample}  
with Exact Hessian Option in Line 9 for maximizing DR-Submodular functions. 
Further, suppose the conditions in 
Assumptions~\ref{assum:constraint}-\ref{assum:second_order} hold, and further 
assume that $F$ is monotone and continuous DR-submodular on 
$\constraint$. Then, if we set the algorithm parameters as $\rho_t = 
(t-1)^{-1}$ and $ \eta_t = T^{-1}$, then the output is a feasible point 
$\bx_{T+1} \in 
\constraint$ and satisfies
\begin{equation*}
\expect[F(\bx_{T+1})] \geq (1-e^{-1})F(\bx^*) - \frac{4R\sqrt{C}}{T^{1/2}} - 
\frac{\bar{L}R^2}{2T},
\end{equation*}
where $C=\max\{4(2G+R\bar{L})^2, 256, [2R(\bar{L}+L)]^4 \}$.
\end{theorem}

Finally, we note that \cref{alg:one_sample} can also be used to solve 
stochastic discrete submodular maximization. Precisely, we can apply 
\cref{alg:one_sample} on the multilinear extension of the discrete submodular 
functions, and round the output to a feasible set by lossless rounding schemes 
like pipage rounding \citep{calinescu2011maximizing} and contention resolution 
method \citep{chekuri2014submodular}.

\section{Oblivious Setting}\label{sec:obl}

In this section, we specifically study the oblivious problem introduced in 
\eqref{eq:ob_op} which is a special case of the non-oblivious problem defined 
in \eqref{eq:formation}. In particular, we show that the proposed one sample 
Frank-Wolfe method becomes significantly simpler under the oblivious setting. 
Also, we show that  the theoretical results for one sample SFW hold under less 
strict assumptions when we are in the oblivious regime. 

\subsection{Algorithm}

As we discussed in Section \ref{sec:sfw}, a major challenge that we face for 
designing a variance reduced Frank-Wolfe method for the non-oblivious setting 
is computing an unbiased estimator of the gradient variation $\Delta_t=\nabla 
F(\bx_t)-\nabla F(\bx_{t-1})$. This is indeed not problematic in the oblivious 
setting, as in this case
$\bz \sim p(\bz)$ is independent of $\bx$ and therefore $\nabla \tF(\bx_t;\bz) 
- \nabla 
\tF(\bx_{t-1};\bz)$ is an unbiased estimator of  the gradient variation 
$\Delta_t = \nabla F(\bx_t) - \nabla F(\bx_{t-1})$. Hence, in the oblivious 
setting, our proposed one sample FW uses the following gradient approximation 
\begin{equation*}\label{eq:momentum_new_2}
\bd_t = (1-\rho_t) (\bd_{t-1}+\tilde{\Delta}_t) + \rho_t \nabla \tF(\bx_t; 
\bz_t),
\end{equation*}
 where $\tilde{\Delta}_t$ is given by 
 $$\tilde{\Delta}_t = \nabla \tF(\bx_t;\bz_t) - \nabla 
		\tF(\bx_{t-1};\bz_t).$$	
 The rest of the algorithm for updating the variable $\bx_t$ is identical to 
 the one for the non-oblivious setting. The description of our proposed 
 algorithm for the oblivious setting is outlined in \cref{alg:one_sample_ob}.

\begin{algorithm}[t]
	\caption{One-Sample SFW (Oblivious Setting)}
	\begin{algorithmic}[1]
		\Require Step sizes $\rho_t \in 
		(0,1),   
		\eta_t \in (0,1)$, initial point $\bx_1 \in 
		\constraint$, total number 
		of iterations $T$
		\Ensure $\bx_{T+1}$ or $\bx_o$, where $\bx_o$ is chosen from $\{\bx_1, 
		\bx_2,\cdots, \bx_T\}$ uniformly at random  
		\For{$t=1,2,\dots, T$}
		\State Sample a point $\bz_t$ according 
		to $p(\bz)$
		\If{$t=1$} 
		\State  Compute $\bd_1 = \nabla \tilde{F}(\bx_1; 
		\bz_1)$		
		\Else
		\State $\tilde{\Delta}_t = \nabla \tF(\bx_t;\bz_t) - \nabla 
		\tF(\bx_{t-1};\bz_t)$
		\State $\bd_t = (1-\rho_t)(\bd_{t-1}+\tilde{\Delta}_t) + \rho_t \nabla 
		\tilde{F}(\bx_t,\bz_t)$
		\EndIf
		\State (non-)convex min.: Update $\bx_{t+1}$ based on 
		\eqref{eq:option1} 
		\State  DR-sub. max.: Update $\bx_{t+1}$ based on \eqref{eq:option2} 
		\EndFor
	\end{algorithmic}
	\label{alg:one_sample_ob}
\end{algorithm}

\vspace{2mm}
\begin{remark}
We note that by rewriting our proposed \sfw method for the
oblivious setting, we recover the variance reduction technique 
applied in the STORM method proposed by \cite{cutkosky2019momentum} with 
different 
settings of parameters. In \citep{cutkosky2019momentum}, however, the STORM 
algorithm was only combined with SGD to solve \emph{unconstrained} non-convex 
minimization problems, while our proposed \sfw method solves convex 
minimization, non-convex minimization, and DR-submodular maximization in a 
\emph{constrained} setting.\end{remark}

\subsection{Theoretical results}

In this section, we show that the variant of one sample stochastic FW for the 
oblivious setting (described in  \cref{alg:one_sample_ob}) recovers the 
theoretical results for the non-oblivious setting with less assumptions. In 
particular, we only require the following condition for the stochastic 
functions $\tilde{F}$ to prove our main results. 

\vspace{2mm}

\begin{assump}\label{assum:function_ob}
The function $\tF$ has uniformly bound gradients, i.e., $ \forall \bx \in 
\constraint, \forall \bz \in 
\mathcal{Z}$
$$\|\nabla \tF(\bx;\bz) \| \leq G.$$
Moreover, the function $\tF$  is uniformly  $L$-smooth, i.e., $\forall \bx, \by 
\in \constraint, \forall \bz \in 
\mathcal{Z}$
$$\|\nabla \tF(\bx;\bz) - \nabla \tF(\by;\bz) \|\leq L \| 
\bx-\by\|$$ 
\end{assump}
\vspace{2mm}

We note that as direct corollaries of \cref{thm:convex,thm:nonconvex,thm:sub}, 
\cref{alg:one_sample_ob} achieves the same optimal convergence rates, which is 
stated in 
\cref{thm:no} formally.

\begin{theorem}\label{thm:no}
Consider the oblivious variant of \sfw outlined in 
Algorithm~\ref{alg:one_sample_ob}, and assume that the conditions in 
\cref{assum:constraint,assum:stoch_bound,assum:function_ob} hold. Then we have
\begin{enumerate}
\item If $F$ is convex on $\constraint$, and we set $\rho_t = 
(t-1)^{-1}$ and $\eta_t = t^{-1}$, then the output is feasible $\bx_{T+1} \in 
\constraint$ and satisfies
\begin{equation*}
\expect[F(\bx_{T+1}) - F(\bx^*)] \leq \mathcal{O}(T^{-1/2}).
\end{equation*}

\item If $F$ is non-convex, and we set $\rho_t 
= (t-1)^{-2/3},$ and $\eta_t = T^{-2/3}$, then the 
output is feasible $\bx_{o} \in \constraint$ and satisfies
\begin{equation*}
\expect[\mathcal{G}(\bx_o)] \leq \mathcal{O}(T^{-1/3}).
\end{equation*}

\item If $F$ is monotone DR-submodular on $\constraint$, and we set $\rho_t = 
(t-1)^{-1}$ and $ \eta_t = T^{-1}$, then the output is a feasible point 
$\bx_{T+1} \in 
\constraint$ and satisfies
\begin{equation*}
\expect[F(\bx_{T+1})] \geq (1-e^{-1})F(\bx^*) - \mathcal{O}(T^{-1/2}).
\end{equation*}
\end{enumerate}
\end{theorem} 

\cref{thm:no} shows that the oblivious version of \sfw requires at most 
$\mathcal{O}(1/\epsilon^2)$ stochastic oracle queries to find an 
$\epsilon$-suboptimal solution for convex minimization, at most 
$\mathcal{O}(1/\epsilon^2)$ stochastic gradient evaluations to achieve a 
$(1-1/e)-\epsilon$ approximate solution for monotone DR-submodular 
maximization, and at most 
$\mathcal{O}(1/\epsilon^3)$ stochastic oracle queries  to find an 
$\epsilon$-first-order stationary point for nonconvex minimization.

\section{Conclusion}
In this paper, we studied the problem of solving constrained  stochastic 
optimization programs  using projection-free methods. We proposed the first 
stochastic variant of the Frank-Wolfe method, called \sfw, that  requires only 
one stochastic sample per iteration while achieving the optimal (or best known) 
complexity bounds for (non-)convex minimization and monotone DR-submodular 
maximization. In particular, we proved that \sfw achieves the 
optimal oracle complexity of $\mathcal{O}(1/\epsilon^2)$ for reaching an  
$\epsilon$-suboptimal solution in the stochastic convex setting, and a 
$(1-1/e)-\epsilon$ approximate solution for a stochastic monotone DR-submodular 
maximization problem. Moreover, in a non-convex setting, \sfw finds an 
$\epsilon$-first-order stationary point after at most 
$\mathcal{O}(1/\epsilon^3)$ iterations, achieving the best known overall 
complexity.  

\appendix
\onecolumn

\section{Proof of \cref{lem:graident_error}}

\begin{proof}
Let $A_t = \|\nabla F(\bx_t) - \bd_t\|^2$. By definition, we have
\begin{equation*}
A_t = \| \nabla F(\bx_{t-1}) - \bd_{t-1}  + \nabla F(\bx_{t}) - 
\nabla F(\bx_{t-1}) - (\bd_t - \bd_{t-1})\|^2.
\end{equation*}
Note that
\begin{equation*}
\bd_t - \bd_{t-1} = - \rho_t \bd_{t-1} + \rho_t \nabla \tF(\bx_t,\bz_t) + 
(1-\rho_t) \tilde{\Delta}_t,
\end{equation*}
and define $\Delta_t = \nabla F(\bx_{t}) - 
\nabla F(\bx_{t-1})$, we have
\begin{equation*}
\begin{split}
A_t &= \|\nabla F(\bx_{t-1}) - \bd_{t-1} + \Delta_t - (1-\rho_t) 
\tilde{\Delta}_t - \rho_t \nabla \tF(\bx_t,\bz_t) + \rho_t \bd_{t-1} \|^2  \\
& = \|\nabla F(\bx_{t-1}) - \bd_{t-1} + (1-\rho_t)(\Delta_t -  
\tilde{\Delta}_t) + \rho_t(\nabla F(\bx_t)-\nabla \tF(\bx_t,\bz_t) + 
\rho_t(\bd_{t-1}-\nabla F(\bx_{t-1}))) \|^2 \\
&= \|(1-\rho_t)(\nabla F(\bx_{t-1}) - \bd_{t-1}) + (1-\rho_t)(\Delta_t -  
\tilde{\Delta}_t) + \rho_t(\nabla F(\bx_t)-\nabla \tF(\bx_t,\bz_t))\|^2
\end{split}
\end{equation*}
Since $\tilde{\Delta}_t$ is an unbiased estimator of $\Delta_t$, $\expect[A_t]$ 
can be decomposed as
\begin{equation}
\label{eq:decomp2}
\begin{split}
\expect[A_t] &= \expect\{(1-\rho_t)^2 \|\nabla F(\bx_{t-1}) - \bd_{t-1} \|^2 + 
(1-\rho_t)^2 \| \Delta_t -  
\tilde{\Delta}_t \|^2  + \rho_t^2 \|\nabla F(\bx_t)-\nabla \tF(\bx_t,\bz_t) 
\|^2 \\
&\quad +2\rho_t(1-\rho_t)\langle\nabla F(\bx_{t-1}) - \bd_{t-1} , \nabla 
F(\bx_t)-\nabla \tF(\bx_t,\bz_t) \rangle \\
&\quad + 2\rho_t(1-\rho_t)\langle \Delta_t -  
\tilde{\Delta}_t , \nabla F(\bx_t)-\nabla \tF(\bx_t,\bz_t) \rangle\}.
\end{split}
\end{equation}
Then we turn to upper bound the items above. First, by \cref{lem:scg++8.1}, we 
have
\begin{equation}
\label{eq:aux1}
\begin{split}
\expect[\|\tilde{\Delta}_t - \Delta_t \|^2] &= \expect[\| 
\tnabla_t^2(\bx_t-\bx_{t-1}) - (\nabla F(\bx_t) - \nabla F(\bx_{t-1}))]\|^2] \\
&\leq \expect[\|\tnabla_t^2(\bx_t-\bx_{t-1}) \|^2] \\
&= \expect[\|\tnabla_t^2(\eta_{t-1}(\bv_{t-1}-\bx_{t-1})) \|^2] \\
&\leq \eta_{t-1}^2D^2 \expect[\| \tnabla_t^2\|^2] \\
&\leq \eta_{t-1}^2D^2\bar{L}^2.
\end{split}
\end{equation}
By Jensen's inequality, we have
\begin{equation}\label{eq:aux2}
\expect[\|\tilde{\Delta}_t - \Delta_t \|] \leq 
\sqrt{\expect[\|\tilde{\Delta}_t - \Delta_t \|^2]} \leq \eta_{t-1} D 
\bar{L},
\end{equation}
and
\begin{equation}\label{eq:aux5}
\expect[\|\nabla F(\bx_t) -\bd_t \|] = \sqrt{\expect[\|\nabla F(\bx_t) -\bd_t 
	\|^2]} = \sqrt{\expect[A_t]}.
\end{equation}

Note that $\bz_t$ is sampled according to $p(\bz; \bx_t(a))$, where $\bx_{t}(a) 
= a\bx_t + (1-a)\bx_{t-1}$. Thus $\nabla \tF(\bx_t,\bz_t)$ is NOT an unbiased 
estimator of $\nabla F(\bx_t)$ when $a \neq 1$, which occurs with probability 
1. 
However, we will show that $\nabla \tF(\bx_t, \bz_t)$ is still a good estimator.
Let $\mathcal{F}_{t-1}$ be the $\sigma$-field generated by all the randomness 
before round $t$, then by Law of Total Expectation, we have
\begin{equation}\label{eq:aux3}
\begin{split}
&\expect[2\rho_t(1-\rho_t)\langle\nabla F(\bx_{t-1}) - \bd_{t-1} , \nabla 
F(\bx_t)-\nabla \tF(\bx_t,\bz_t) \rangle] \\
=& \expect[\expect[2\rho_t(1-\rho_t)\langle\nabla F(\bx_{t-1}) - \bd_{t-1} , 
\nabla 
F(\bx_t)-\nabla \tF(\bx_t,\bz_t) \rangle|\mathcal{F}_{t-1},\bx_t(a)]] \\
=& \expect[2\rho_t(1-\rho_t)\langle\nabla F(\bx_{t-1}) - \bd_{t-1} , 
\expect[\nabla 
F(\bx_t)-\nabla \tF(\bx_t,\bz_t) |\mathcal{F}_{t-1},\bx_t(a)]\rangle],
\end{split}
\end{equation}
where 
\begin{equation*}
\expect[\nabla 
F(\bx_t)-\nabla \tF(\bx_t,\bz_t) |\mathcal{F}_{t-1}]\rangle] = \nabla F(\bx_t) 
- \nabla F(\bx_t(a)) + \nabla F(\bx_t(a)) - \expect[\nabla 
\tF(\bx_t,\bz_t)|\mathcal{F}_{t-1},\bx_t(a)].
\end{equation*}
By \cref{lem:scg++8.1}, $F$ is $\bar{L}$-smooth, thus 
\begin{equation*}
\|\nabla F(\bx_t) - \nabla F(\bx_t(a)) \| \leq \bar{L} \|\bx_t - \bx_t(a) \| = 
\bar{L}(1-a)\| \eta_{t-1}(\bv_{t-1}-\bx_{t-1})\| \leq \eta_{t-1}D\bar{L}.
\end{equation*}
We also have
\begin{equation*}
\begin{split}
\| \nabla F(\bx_t(a)) - \expect[\nabla 
\tF(\bx_t,\bz_t)|\mathcal{F}_{t-1},\bx_t(a)]\| & = \| \int [\nabla \tF 
(\bx_t(a); \bz) - 
\nabla \tF (\bx_t; \bz)] p(\bz; \bx_t(a)) \mathrm{d}\bz \|\\
&\leq \int \| \nabla \tF (\bx_t(a); \bz) - 
\nabla \tF (\bx_t; \bz)\| p(\bz; \bx_t(a)) \mathrm{d}\bz \\
&\leq \int L \|\bx_t(a)-\bx_t \| p(\bz; \bx_t(a)) \mathrm{d}\bz \\
&\leq \eta_{t-1}DL,
\end{split}
\end{equation*}
where the second inequality holds because of \cref{assum:second_order}. Combine 
the analysis above with \cref{eq:aux3},
we have
\begin{equation}\label{eq:aux7}
\begin{split}
&\expect[2\rho_t(1-\rho_t)\langle\nabla F(\bx_{t-1}) - \bd_{t-1} , \nabla 
F(\bx_t)-\nabla \tF(\bx_t,\bz_t) \rangle] \\
\leq & \expect[2\rho_t(1-\rho_t)\|\nabla F(\bx_{t-1}) - \bd_{t-1}\| \cdot 
\|\expect[\nabla 
F(\bx_t)-\nabla \tF(\bx_t,\bz_t) |\mathcal{F}_{t-1}]\|] \\
\leq & 2\rho_t(1-\rho_t) \expect[\|\nabla F(\bx_{t-1}) - \bd_{t-1}\|] \cdot 
(\eta_{t-1}D\bar{L} + \eta_{t-1}DL) \\
\leq & 2\eta_{t-1}\rho_t(1-\rho_t) \sqrt{\expect[A_{t-1}]}D(\bar{L}+L).
\end{split}
\end{equation}

Finally, by \cref{assum:gradient_norm}, we have $\| \nabla F(\bx_t)-\nabla 
\tF(\bx_t,\bz_t) \| \leq 2G$. Thus 
\begin{equation}\label{eq:aux6}
\rho_t^2 \|\nabla F(\bx_t)-\nabla \tF(\bx_t,\bz_t) 
\|^2 \leq 4\rho_t^2G^2,
\end{equation}
and
\begin{equation}\label{eq:aux4}
\begin{split}
\expect [2\rho_t(1-\rho_t)\langle \Delta_t -  
\tilde{\Delta}_t , \nabla F(\bx_t)-\nabla \tF(\bx_t,\bz_t) \rangle] &\le 
\expect[2\rho_t(1-\rho_t)\| \Delta_t -  
\tilde{\Delta}_t\| \cdot \|\nabla F(\bx_t)-\nabla \tF(\bx_t,\bz_t)] \| \\
&\leq 4\eta_{t-1}\rho_t(1-\rho_t)GD\bar{L}.
\end{split}
\end{equation}
Combine \cref{eq:decomp2,eq:aux1,eq:aux6,eq:aux7,eq:aux4}, we have
\begin{equation*}
\begin{split}
\expect[A_t] &\leq (1-\rho_t)^2 \expect[A_{t-1}] + (1-\rho_t)^2 
\eta_{t-1}^2D^2\bar{L}^2 + \rho_t^2 4G^2 + 
2\eta_{t-1}\rho_t(1-\rho_t) \sqrt{\expect[A_{t-1}]}D(\bar{L}+L) \\
&\quad + 4\eta_{t-1}\rho_t(1-\rho_t)GD\bar{L}.
\end{split}
\end{equation*}
For the simplicity of analysis, we replace $t$ by $t+1$, and have
\begin{equation}\label{eq:recursive}
\begin{split}
&\expect[A_{t+1}] \\
\leq& (1-\rho_{t+1})^2 \expect[A_{t}] + (1-\rho_{t+1})^2 
\eta_{t}^2D^2\bar{L}^2 + \rho_{t+1}^2 4G^2 + 
2\eta_{t}\rho_{t+1}(1-\rho_{t+1}) \sqrt{\expect[A_{t}]}D(\bar{L}+L)\\
&\quad + 4\eta_{t}\rho_{t+1}(1-\rho_{t+1})GD\bar{L} \\
\leq& (1-\frac{1}{t^\alpha})^2 \expect[A_{t}] + 
\frac{D^2\bar{L}^2+4G^2+4GD\bar{L}}{t^{2\alpha}} + 
\frac{2D(\bar{L}+L)}{t^{2\alpha}}\sqrt{\expect[A_{t}]}.
\end{split}
\end{equation}
We claim that $\expect[A_{t}] \leq C t^{-\alpha}$, and prove it by induction. 
Before the proof, we first analyze one item in the definition of $C: 
\frac{2(2G+D\bar{L})^2}{2-2^{-\alpha}-\alpha}$. Define $h(\alpha) = 
2-2^{-\alpha}-\alpha$. Since $h'(\alpha) = 2^{-\alpha}\ln(2)-1\le 0$ for 
$\alpha \in (0,1]$, so $1=h(0) \ge h(\alpha) \ge h(1) = 1/2 >0, \forall \alpha 
\in (0,1]$. As a result, $2 \le \frac{2}{2-2^{-\alpha}-\alpha}\le 4$.

When $t=1$, we have
\begin{equation*}
\expect[A_{1}] = \expect[\|\nabla F(\bx_1)-\nabla \tF (\bx_1;\bz_1) \|^2] \leq 
(2G)^2 \leq \frac{2(2G+D\bar{L})^2}{2-2^{-\alpha}-\alpha} / 1 \leq C\cdot 
1^{-\alpha}
\end{equation*}
When $t=2$, since $\rho_2 = 1$, we have
\begin{equation*}
\begin{split}
\expect[A_{2}] = \expect[\|\nabla 
\tF(\bx_2,\bz_2) - \nabla F(\bx_2)\|^2] \leq (2G)^2 \leq 
\frac{2(2G+D\bar{L})^2}{2-2^{-\alpha}-\alpha} / 2 \leq C \cdot 2^{-\alpha}.
\end{split}
\end{equation*}
Now assume for $t \ge 2$, we have $\expect[A_{t}] \leq C t^{-\alpha}$, by 
\cref{eq:recursive} and the definition of $C$, we have
\begin{equation}\label{eq:induction}
\begin{split}
\expect[A_{t+1}] &\leq (1-\frac{1}{t^\alpha})^2 \cdot Ct^{-\alpha} + 
\frac{(2G+D\bar{L})^2}{t^{2\alpha}} + 
\frac{2D(\bar{L}+L)}{t^{(5/2)\alpha}}\sqrt{C}\\
&\leq Ct^{-\alpha} - 2Ct^{-2\alpha} + Ct^{-3\alpha} + 
\frac{(2-2^{-\alpha}-\alpha)C}{2t^{2\alpha}} + \frac{C^{3/4}}{t^{(5/2)\alpha}} 
\\ 
&\leq \frac{C}{t^\alpha} + \frac{-2C+Ct^{-\alpha}+ 
(2-2^{-\alpha}-\alpha)C/2+t^{-\alpha/2}C 
/ C^{1/4}}{t^{2\alpha}} \\
&\leq \frac{C}{t^\alpha} + \frac{C[-2+2^{-\alpha} +(2-2^{-\alpha}-\alpha)/2 
+(2-2^{-\alpha}-\alpha)/2]}{t^{2\alpha}} \\
&\leq \frac{C}{t^\alpha} - \frac{\alpha C}{t^{2\alpha}}.
\end{split}
\end{equation}
Define $g(t) = t^{-\alpha}$, then $g(t)$ is a convex function for $\alpha \in 
(0,1]$. Thus we have $g(t+1) - g(t) \ge g'(t)$, \emph{i.e.}, $(t+1)^{-\alpha} 
- t^{-\alpha} \geq -\alpha t^{-(\alpha+1)}$. So we have
\begin{equation}
\frac{C}{t^{\alpha}}- \frac{\alpha C}{t^{2\alpha}} \leq C(t^{-\alpha} - \alpha 
t^{-(1+\alpha)}) \leq C (t+1)^{-\alpha}.
\end{equation}
Combine with \cref{eq:induction}, we have $\expect[A_{t+1}] \leq C 
(t+1)^{-\alpha}$. Thus by induction, we have $\expect[A_{t}] \leq C 
t^{-\alpha}, \forall t \ge 1$.
\end{proof}

\section{Proof of \cref{lem:graident_error_2}}
The only difference with the proof of \cref{lem:graident_error} is the bound 
for $\expect \| \tilde{\Delta}_t - \Delta_t\|$. Specifically, we have
\begin{equation*}
\begin{split}
\expect[\|\tilde{\Delta}_t - \Delta_t \|^2] &= \expect[\| 
\tilde{\Delta}_t - \tnabla_t^2(\bx_t-\bx_{t-1}) + \tnabla_t^2(\bx_t-\bx_{t-1}) 
- (\nabla F(\bx_t) - \nabla F(\bx_{t-1}))]\|^2] \\
&= \expect[\|\tilde{\Delta}_t - \tnabla_t^2(\bx_t-\bx_{t-1}) \|^2] + 
\expect[\|\tnabla_t^2(\bx_t-\bx_{t-1}) 
- (\nabla F(\bx_t) - \nabla F(\bx_{t-1})) \|^2]\\
&\leq [D^2L_2\delta_t(1+\tF(\bx_t(a),\bz_t))]^2 + \eta_{t-1}^2D^2\bar{L}^2\\
&\leq (1+B)^2L_2^2D^4\delta_t^2 + \eta_{t-1}^2D^2\bar{L}^2 \\
&\leq 4\eta_{t-1}^2D^2\bar{L}^2.
\end{split}
\end{equation*}

Then by the analysis same to the proof of \cref{lem:graident_error}, we have
\begin{equation*}
\expect[A_{t+1}]\leq  (1-\frac{1}{t^\alpha})^2 \expect[A_{t}] + 
\frac{4(D^2\bar{L}^2+G^2+GD\bar{L})}{t^{2\alpha}} + 
\frac{4D(\bar{L}+L)}{t^{2\alpha}}\sqrt{\expect[A_{t}]},
\end{equation*}
and thus $\expect[A_{t+1}] \leq C 
(t+1)^{-\alpha}$, where $C\! =\! \max \! \left\{ 
\!\frac{8(D^2\bar{L}^2+G^2+GD\bar{L})}{2\!-\!2^{-\alpha}\!-\!\alpha},
\left[\frac{2}{2\!-\!2^{-\alpha}\!-\!\alpha} \right]^4\!\!, 
[4D(\bar{L}\!+\!L)]^4 \! \right\}\!.$

\section{Proof of \cref{thm:convex}}
First, since $\bx_{t+1} = (1-\eta_t)\bx_t + \eta_t \bv_t$ is a convex 
combination of $\bx_t, \bv_t$, and $\bx_1 \in \constraint, \bv_t \in 
\constraint, \forall\ t$, we can prove $\bx_t \in \constraint, \forall\ t $ by 
induction. So $\bx_{T+1} \in \constraint$.

Then we present an auxiliary lemma.
\begin{lemma}\label{lem:aryan1}
Under the condition of \cref{thm:convex}, in \cref{alg:one_sample}, we have 
\begin{equation*}
F(\bx_{t+1}) - F(\bx^*) \leq (1-\eta_t)(F(\bx_t) - F(\bx^*)) + \eta_tD\|\nabla 
F(\bx_t) - \bd_t\| + \frac{\bar{L}D^2\eta_t^2}{2}.
\end{equation*}
\end{lemma}

By Jensen's inequality and \cref{lem:graident_error} with $\alpha=1$, we have
\begin{equation*}
\expect[\| \nabla F(\bx_t) - \bd_t\| ] \leq \sqrt{\expect[\| \nabla F(\bx_t) - 
\bd_t\|^2 ]} \leq \frac{\sqrt{C}}{\sqrt{t}},
\end{equation*}
where $C=\max \{4(2G+D\bar{L})^2, 256, [2D(\bar{L}+L)]^4 \}$. Then by 
\cref{lem:aryan1}, we have
\begin{equation}\label{eq:convex1}
\begin{split}
&\expect[F(\bx_{T+1}) - F(\bx^*)] \\
\leq& (1-\eta_T)\expect[F(\bx_T) - F(\bx^*)] + 
\eta_TD\expect[\|\nabla 
F(\bx_T) - \bd_T\|] + \frac{\bar{L}D^2\eta_T^2}{2} \\
={}& \prod_{i=1}^T (1-\eta_i) \expect[F(\bx_1) - F(\bx^*)] + D\sum_{k=1}^T 
\eta_k\expect[\|\nabla F(\bx_k) - \bd_k\|] \prod_{i=k+1}^{T} (1-\eta_i)\\
&\quad + 
\frac{\bar{L}D^2}{2} \sum_{k=1}^T \eta_k^2\prod_{i=k+1}^{T} (1-\eta_i) \\
\leq{}& 0 + D\sum_{k=1}^T k^{-1}\frac{\sqrt{C}}{\sqrt{k}} \prod_{i=k+1}^{T} 
\frac{i-1}{i} + \frac{\bar{L}D^2}{2} \sum_{k=1}^T k^{-2} 
\prod_{i=k+1}^{T}\frac{i-1}{i} \\
={}& \frac{\sqrt{C}D}{T}\sum_{k=1}^T \frac{1}{\sqrt{k}} + 
\frac{\bar{L}D^2}{2T}\sum_{k=1}^T k^{-1}.
\end{split}
\end{equation}
Since
\begin{equation*}
\sum_{k=1}^T \frac{1}{\sqrt{k}} \leq \int_{0}^T x^{-1/2}\mathrm{d}x = 2\sqrt{T},
\end{equation*}
and
\begin{equation*}
\sum_{k=1}^T k^{-1} \leq 1 + \int_1^T x^{-1}\mathrm{d}x = 1+\ln T,
\end{equation*}
by \cref{eq:convex1}, we have
\begin{equation*}
\expect[F(\bx_{T+1}) - F(\bx^*)] \leq \frac{2\sqrt{C}D}{\sqrt{T}} + 
\frac{\bar{L}D^2}{2T}(1+\ln T).
\end{equation*}

\section{Proof of \cref{thm:nonconvex}}
First, since $\bx_{t+1} = (1-\eta_t)\bx_t + \eta_t \bv_t$ is a convex 
combination of $\bx_t, \bv_t$, and $\bx_1 \in \constraint, v_t \in \constraint, 
\forall\ t$, we can prove $\bx_t \in \constraint, \forall\ t $ by induction. So 
$\bx_{o} \in \constraint$.

Note that if we define $\bv_t^\prime = \argmin_{\bv \in \constraint}\langle 
\bv, 
\nabla F(\bx_t)\rangle$, then $\mathcal{G}(\bx_t)=\langle 
\bv^\prime_t-\bx_t,-\nabla F(\bx_t)\rangle = -\langle \bv^\prime_t-\bx_t,\nabla 
F(\bx_t)\rangle$. So we have

\begin{equation*}
\begin{split}
F(\bx_{t+1}) 
\stackrel{(a)}{\leq}& F(\bx_t) + \langle \nabla f(\bx_t),\bx_{t+1}-\bx_t 
\rangle + 
\frac{\bar{L}}{2}\|\bx_{t+1}-\bx_t\|^2 \\
={}&F(\bx_t) + \langle \nabla F(\bx_t),\eta_t(\bv_t-\bx_t) \rangle + 
\frac{\bar{L}}{2}\|\eta_t(\bv_t-\bx_t)\|^2 \\
\stackrel{(b)}{\leq}& F(\bx_t) + \eta_t \langle \nabla F(\bx_t),\bv_t-\bx_t 
\rangle+\frac{\bar{L}\eta_t^2D^2}{2} \\
={}& F(\bx_t) + \eta_t \langle  \bd_t,\bv_t-\bx_t \rangle+ \eta_t \langle 
\nabla 
F(\bx_t)-\bd_t,\bv_t-\bx_t \rangle  + \frac{\bar{L}\eta_t^2D^2}{2} \\
\stackrel{(c)}{\leq}& F(\bx_t) + \eta_t \langle  \bd_t,\bv^\prime_t-\bx_t 
\rangle+ \eta_t \langle \nabla F(\bx_t)-\bd_t,\bv_t-\bx_t \rangle + 
\frac{\bar{L}\eta_t^2D^2}{2} \\
={}& F(\bx_t) + \eta_t \langle \nabla F(\bx_t),\bv^\prime_t-\bx_t \rangle  + 
\eta_t 
\langle  \bd_t-\nabla F(\bx_t),\bv^\prime_t-\bx_t \rangle \\
&\quad + \eta_t \langle 
\nabla F(\bx_t)-\bd_t,\bv_t-\bx_t \rangle + 
\frac{\bar{L}\eta_t^2D^2}{2} \\
={}& F(\bx_t) - \eta_t \mathcal{G}(\bx_t) + \eta_t \langle \nabla 
F(\bx_t)-\bd_t,\bv_t-\bv^\prime_t \rangle + \frac{\bar{L}\eta_t^2D^2}{2} \\
\stackrel{(d)}{\leq}& F(\bx_t) - \eta_t \mathcal{G}(\bx_t) + \eta_t \|\nabla 
F(\bx_t)-\bd_t\|\| \bv_t-\bv^\prime_t\| + \frac{\bar{L}\eta_t^2D^2}{2} \\
\stackrel{(e)}{\leq}& F(\bx_t) - \eta_t \mathcal{G}(\bx_t) + \eta_tD \|\nabla 
F(\bx_t)-\bd_t\| + \frac{\bar{L}\eta_t^2D^2}{2},
\end{split}
\end{equation*}
where we used the fact that $F$ is $\bar{L}$-smooth in inequality (a). 
Inequalities 
(b), (e) hold because of \cref{assum:constraint}. Inequality (c) is due to the 
optimality of $\bv_t$, and 
in (d), we applied the Cauchy-Schwarz inequality.

Rearrange the inequality above, we have 
\begin{equation}
\label{eq:bound_on_individual_gap}
\eta_t \mathcal{G}(\bx_t) \leq F(\bx_t)- F(\bx_{t+1})+ \eta_tD \|\nabla 
F(\bx_t)-\bd_t\| + \frac{\bar{L}\eta_t^2D^2}{2}.   
\end{equation}

Apply \cref{eq:bound_on_individual_gap} recursively for $t = 1, 2, \cdots, 
T$, and take expectations, we attain the following inequality:
\begin{equation*}
\begin{split}
\sum_{t=1}^T \eta_t \expect[\mathcal{G}(\bx_t)]
\leq F(\bx_1)-F(\bx_{T+1})
+D\sum_{t=1}^T \eta_t\expect[\|\nabla F(\bx_t)-\bd_t\|] 
+\frac{\bar{L}D^2}{2}\sum_{t=1}^T \eta_t^2.
\end{split}
\end{equation*}
By Jensen's inequality \cref{lem:graident_error} with $\alpha=2/3$, we have
\begin{equation*}
\expect[\|\nabla F(\bx_t)-\bd_t\|] \leq \sqrt{\expect[\|\nabla 
F(\bx_t)-\bd_t\|^2]} \leq \frac{\sqrt{C}}{t^{1/3}},
\end{equation*}
where $C = \max \{\frac{2(2G+D\bar{L})^2}{4/3-2^{-2/3}}, \left( 
\frac{2}{4/3-2^{-2/3}}\right)^4, [2D(\bar{L}+L)]^4 \}$. Since $\eta_t = 
T^{-2/3}$, we have
\begin{equation*}
\begin{split}
\expect[\mathcal{G}(\bx_o)] &= \frac{\sum_{t=1}^T 
\expect[\mathcal{G}(\bx_t)]}{T} \\
&\leq \frac{1}{T\cdot T^{-2/3}} [F(\bx_1) - F(\bx_{T+1}) + D\sum_{t=1}^T 
T^{-2/3}\frac{\sqrt{C}}{t^{1/3}} + \frac{\bar{L}D^2}{2}\sum_{t=1}^T T^{-4/3}] \\
&\leq \frac{1}{T^{1/3}}[2B + 
D\sqrt{C}T^{-2/3}\frac{3}{2}T^{2/3}+\frac{\bar{L}D^2}{2T^{1/3}}] \\
&= \frac{2B+3\sqrt{C}D/2}{T^{1/3}} + \frac{\bar{L}D^2}{2T^{2/3}}, 
\end{split}
\end{equation*}
where the second inequality holds because $\sum_{t=1}^T t^{-1/3} \leq 
\int_{0}^T x^{-1/3}\mathrm{d}x = \frac{3}{2}T^{2/3}$.

\section{Proof of \cref{thm:sub}}
First, since $\bx_{t+1} = \bx_t + \eta_t \bv_t = \bx_t + T^{-1} \bv_t$, we have
$\bx_{T+1} = \frac{\sum_{t=1}^T \bv_t}{T} \in \constraint$. Also, because now 
$\|\bx_{t+1} -\bx_{t} \| = \|\eta_t \bv_t \| \leq \eta_t R$, (rather than 
$\eta_t D$), \cref{lem:graident_error} holds with new constant $C = \max 
\{\frac{2(2G+R\bar{L})^2}{2-2^{-\alpha}-\alpha}, 
\left(\frac{2}{2-2^{-\alpha}-\alpha} \right)^4, [2R(\bar{L}+L)]^4 \}$. Since 
$\alpha=1$, we have $C=\max\{4(2G+R\bar{L})^2, 256, [2R(\bar{L}+L)]^4 \}$.
Then by Jensen's inequality, we have
\begin{equation*}
\expect[\|\nabla F(\bx_t) -\bd_t \|] \leq \sqrt{\expect[\|\nabla F(\bx_t) 
-\bd_t \|^2]}\leq \frac{\sqrt{C}}{\sqrt{t}}.
\end{equation*}	

We observe that 
\begin{equation}\label{eq:sub_aux1}
\begin{split}
F(\bx_{t+1}) &\stackrel{(a)}{\geq} F(\bx_t) + \langle \nabla F(\bx_t), 
\bx_{t+1}-\bx_t \rangle -\frac{\bar{L}}{2} \|\bx_{t+1} -\bx_t \| \\
&= F(\bx_t) + \frac{1}{T} \langle \nabla F(\bx_t), \bv_t \rangle - 
\frac{\bar{L}}{2T^2} \|\bv_t \| \\
&\stackrel{(b)}{\geq} F(\bx_t) + \frac{1}{T} \langle \bd_t, \bv_t \rangle + 
\frac{1}{T} \langle \nabla F(\bx_t) -\bd_t, \bv_t \rangle - 
\frac{\bar{L}R^2}{2T^2} \\
&\stackrel{(c)}{\geq} F(\bx_t) + \frac{1}{T} \langle \bd_t, \bx^* \rangle + 
\frac{1}{T} \langle \nabla F(\bx_t) -\bd_t, \bv_t \rangle - 
\frac{\bar{L}R^2}{2T^2} \\
&= F(\bx_t) + \frac{1}{T} \langle \nabla F(\bx_t), \bx^* \rangle + 
\frac{1}{T} \langle \nabla F(\bx_t) -\bd_t, \bv_t - \bx^*\rangle - 
\frac{\bar{L}R^2}{2T^2} \\
&\stackrel{(d)}{\geq} F(\bx_t) + \frac{F(\bx^*)-F(\bx_t)}{T}- \frac{1}{T} 
\langle \nabla F(\bx_t) -\bd_t, -\bv_t + \bx^*\rangle - \frac{\bar{L}R^2}{2T^2} 
\\
&\stackrel{(e)}{\geq} F(\bx_t) + \frac{F(\bx^*)-F(\bx_t)}{T}- \frac{1}{T} 
\|\nabla F(\bx_t) -\bd_t\|\cdot \|-\bv_t + \bx^*\| - \frac{\bar{L}R^2}{2T^2} \\
&\stackrel{(f)}{\geq}  F(\bx_t) + \frac{F(\bx^*)-F(\bx_t)}{T}- \frac{1}{T} 
2R\|\nabla F(\bx_t) -\bd_t\| - \frac{\bar{L}R^2}{2T^2},
\end{split} 
\end{equation}
where inequality $(a)$ holds because of the $\bar{L}$-smoothness of $F$, 
inequalities $(b), (e)$ comes from \cref{assum:constraint}. We used the 
optimality of $\bv_t$ in inequality $(c)$, and applied the Cauchy-Schwarz 
inequality in$(e)$. Inequality $(d)$ is a little involved, since $F$ is 
monotone and concave in positive directions, we have
\begin{equation*}
F(\bx^*) - F(\bx_t) \leq F(\bx^* \vee \bx_t) - F(\bx_t) \leq \langle \nabla 
F(\bx_t), \bx^* \vee \bx_t - \bx_t \rangle = \langle \nabla F(\bx_t), (\bx^*- 
\bx_t) \vee 0 \rangle \leq \langle \nabla F(\bx_t), \bx^* \rangle.
\end{equation*}  

Taking expectations on both sides of \cref{eq:sub_aux1},
\begin{equation*}
\expect[F(\bx_{t+1})] \geq \expect[F(\bx_t)] + 
\frac{F(\bx^*)-\expect[F(\bx_t)]}{T}- \frac{2R}{T}\frac{\sqrt{C}}{\sqrt{t}}  - 
\frac{\bar{L}R^2}{2T^2}.
\end{equation*}
Or 
\begin{equation*}
F(\bx^*) - \expect[F(\bx_{t+1})] \leq 
(1-\frac{1}{T})[F(\bx^*)-\expect[F(\bx_t)]] + 
\frac{2R}{T}\frac{\sqrt{C}}{\sqrt{t}} + \frac{\bar{L}R^2}{2T^2}
\end{equation*}
Apply the inequality above recursively for $t = 1, 2, \cdots, T$, we have
\begin{equation*}
\begin{split}
F(\bx^*) - \expect[F(\bx_{T+1})] &\leq (1-\frac{1}{T})^T [F(\bx^*)-F(\bx_1)] + 
\frac{2R\sqrt{C}}{T} \sum_{t=1}^T t^{-1/2} + \frac{\bar{L}R^2}{2T} \\
&\leq e^{-1} F(\bx^*) + \frac{4R\sqrt{C}}{T^{1/2}} + \frac{\bar{L}R^2}{2T},
\end{split}
\end{equation*}
where the second inequality holds since $\sum_{t=1}^T t^{-1/2} \leq \int_{0}^T 
x^{-1/2}\mathrm{d}x = 2T^{1/2}$. Thus we have
\begin{equation*}
\expect[F(\bx_{T+1})] \geq (1-e^{-1})F(\bx^*) - \frac{4R\sqrt{C}}{T^{1/2}} - 
\frac{\bar{L}R^2}{2T}.
\end{equation*}

\bibliographystyle{plainnat}
\bibliography{sfw,scgpp}

\begin{thebibliography}{44}
\providecommand{\natexlab}[1]{#1}
\providecommand{\url}[1]{\texttt{#1}}
\expandafter\ifx\csname urlstyle\endcsname\relax
  \providecommand{\doi}[1]{doi: #1}\else
  \providecommand{\doi}{doi: \begingroup \urlstyle{rm}\Url}\fi

\bibitem[Allen-Zhu and Hazan(2016)]{allen2016variance}
Zeyuan Allen-Zhu and Elad Hazan.
\newblock Variance reduction for faster non-convex optimization.
\newblock In \emph{International Conference on Machine Learning}, pages
  699--707, 2016.

\bibitem[Bach(2015)]{bach2015submodular}
Francis Bach.
\newblock Submodular functions: from discrete to continous domains.
\newblock \emph{arXiv preprint arXiv:1511.00394}, 2015.

\bibitem[Bellet et~al.(2015)Bellet, Liang, Garakani, Balcan, and
  Sha]{bellet2015distributed}
Aur{\'e}lien Bellet, Yingyu Liang, Alireza~Bagheri Garakani, Maria-Florina
  Balcan, and Fei Sha.
\newblock A distributed frank-wolfe algorithm for communication-efficient
  sparse learning.
\newblock In \emph{Proceedings of the 2015 SIAM International Conference on
  Data Mining}, pages 478--486. SIAM, 2015.

\bibitem[Bian et~al.(2017{\natexlab{a}})Bian, Levy, Krause, and
  Buhmann]{bian2017continuous}
An~Bian, Kfir Levy, Andreas Krause, and Joachim~M Buhmann.
\newblock Continuous dr-submodular maximization: Structure and algorithms.
\newblock In \emph{Advances in Neural Information Processing Systems}, pages
  486--496, 2017{\natexlab{a}}.

\bibitem[Bian et~al.(2017{\natexlab{b}})Bian, Mirzasoleiman, Buhmann, and
  Krause]{bian16guaranteed}
An~Bian, Baharan Mirzasoleiman, Joachim~M. Buhmann, and Andreas Krause.
\newblock Guaranteed non-convex optimization: Submodular maximization over
  continuous domains.
\newblock In \emph{AISTATS}, February 2017{\natexlab{b}}.

\bibitem[Bian et~al.(2018)Bian, Buhmann, and Krause]{bian2018optimal}
An~Bian, Joachim~M Buhmann, and Andreas Krause.
\newblock Optimal dr-submodular maximization and applications to provable mean
  field inference.
\newblock \emph{arXiv preprint arXiv:1805.07482}, 2018.

\bibitem[Braun et~al.(2017)Braun, Pokutta, and Zink]{braun2017lazifying}
G\'{a}bor Braun, Sebastian Pokutta, and Daniel Zink.
\newblock Lazifying conditional gradient algorithms.
\newblock In \emph{Proceedings of the 34th International Conference on Machine
  Learning - Volume 70}, ICML'17, pages 566--575, 2017.

\bibitem[Calinescu et~al.(2011)Calinescu, Chekuri, P{\'a}l, and
  Vondr{\'a}k]{calinescu2011maximizing}
Gruia Calinescu, Chandra Chekuri, Martin P{\'a}l, and Jan Vondr{\'a}k.
\newblock Maximizing a monotone submodular function subject to a matroid
  constraint.
\newblock \emph{SIAM Journal on Computing}, 40\penalty0 (6):\penalty0
  1740--1766, 2011.

\bibitem[Chekuri et~al.(2014)Chekuri, Vondr{\'a}k, and
  Zenklusen]{chekuri2014submodular}
Chandra Chekuri, Jan Vondr{\'a}k, and Rico Zenklusen.
\newblock Submodular function maximization via the multilinear relaxation and
  contention resolution schemes.
\newblock \emph{SIAM Journal on Computing}, 43\penalty0 (6):\penalty0
  1831--1879, 2014.

\bibitem[Chen et~al.(2018)Chen, Hassani, and Karbasi]{Chen2018Online}
Lin Chen, Hamed Hassani, and Amin Karbasi.
\newblock Online continuous submodular maximization.
\newblock In \emph{AISTATS}, pages 1896--1905, 2018.

\bibitem[Chen et~al.(2019)Chen, Zhang, Hassani, and Karbasi]{chen2019black}
Lin Chen, Mingrui Zhang, Hamed Hassani, and Amin Karbasi.
\newblock Black box submodular maximization: Discrete and continuous settings.
\newblock \emph{arXiv preprint arXiv:1901.09515}, 2019.

\bibitem[Cutkosky and Orabona(2019)]{cutkosky2019momentum}
Ashok Cutkosky and Francesco Orabona.
\newblock Momentum-based variance reduction in non-convex sgd.
\newblock \emph{arXiv preprint arXiv:1905.10018}, 2019.

\bibitem[Defazio and Bottou(2018)]{defazio2018ineffectiveness}
Aaron Defazio and L{\'e}on Bottou.
\newblock On the ineffectiveness of variance reduced optimization for deep
  learning.
\newblock \emph{arXiv preprint arXiv:1812.04529}, 2018.

\bibitem[Fang et~al.(2018)Fang, Li, Lin, and Zhang]{fang2018spider}
Cong Fang, Chris~Junchi Li, Zhouchen Lin, and Tong Zhang.
\newblock Spider: Near-optimal non-convex optimization via stochastic
  path-integrated differential estimator.
\newblock In \emph{Advances in Neural Information Processing Systems}, pages
  687--697, 2018.

\bibitem[Frank and Wolfe(1956)]{frank1956algorithm}
Marguerite Frank and Philip Wolfe.
\newblock An algorithm for quadratic programming.
\newblock \emph{Naval Research Logistics (NRL)}, 3\penalty0 (1-2):\penalty0
  95--110, 1956.

\bibitem[Garber and Hazan(2015)]{garber2015faster}
Dan Garber and Elad Hazan.
\newblock Faster rates for the frank-wolfe method over strongly-convex sets.
\newblock In \emph{ICML}, volume~15, pages 541--549, 2015.

\bibitem[Hassani et~al.(2017)Hassani, Soltanolkotabi, and
  Karbasi]{hassani2017gradient}
Hamed Hassani, Mahdi Soltanolkotabi, and Amin Karbasi.
\newblock Gradient methods for submodular maximization.
\newblock \emph{arXiv preprint arXiv:1708.03949}, 2017.

\bibitem[Hassani et~al.(2019)Hassani, Karbasi, Mokhtari, and
  Shen]{hassani2019stochastic}
Hamed Hassani, Amin Karbasi, Aryan Mokhtari, and Zebang Shen.
\newblock Stochastic conditional gradient++.
\newblock \emph{arXiv preprint arXiv:1902.06992}, 2019.

\bibitem[Hazan and Kale(2012)]{DBLP:conf/icml/HazanK12}
Elad Hazan and Satyen Kale.
\newblock Projection-free online learning.
\newblock In \emph{Proceedings of the 29th International Conference on Machine
  Learning, {ICML} 2012, Edinburgh, Scotland, UK, June 26 - July 1, 2012},
  pages 1843--1850, 2012.

\bibitem[Hazan and Luo(2016)]{hazan2016variance}
Elad Hazan and Haipeng Luo.
\newblock Variance-reduced and projection-free stochastic optimization.
\newblock In \emph{ICML}, pages 1263--1271, 2016.

\bibitem[Jaggi(2013)]{jaggi2013revisiting}
Martin Jaggi.
\newblock Revisiting frank-wolfe: Projection-free sparse convex optimization.
\newblock In \emph{ICML}, pages 427--435, 2013.

\bibitem[Johnson and Zhang(2013)]{Johnson2013Accelerating}
Rie Johnson and Tong Zhang.
\newblock Accelerating stochastic gradient descent using predictive variance
  reduction.
\newblock In \emph{NIPS}, pages 315--323, 2013.

\bibitem[Kulesza et~al.(2012)Kulesza, Taskar, et~al.]{kulesza2012determinantal}
Alex Kulesza, Ben Taskar, et~al.
\newblock Determinantal point processes for machine learning.
\newblock \emph{Foundations and Trends{\textregistered} in Machine Learning},
  5\penalty0 (2--3):\penalty0 123--286, 2012.

\bibitem[Lacoste-Julien(2016)]{lacoste2016convergence}
Simon Lacoste-Julien.
\newblock Convergence rate of frank-wolfe for non-convex objectives.
\newblock \emph{arXiv preprint arXiv:1607.00345}, 2016.

\bibitem[Lacoste-Julien and Jaggi(2015)]{lacoste2015global}
Simon Lacoste-Julien and Martin Jaggi.
\newblock On the global linear convergence of frank-wolfe optimization
  variants.
\newblock In \emph{Advances in Neural Information Processing Systems}, pages
  496--504, 2015.

\bibitem[Lafond et~al.(2016)Lafond, Wai, and Moulines]{lafond2016d}
Jean Lafond, Hoi-To Wai, and Eric Moulines.
\newblock D-fw: Communication efficient distributed algorithms for
  high-dimensional sparse optimization.
\newblock In \emph{Acoustics, Speech and Signal Processing (ICASSP), 2016 IEEE
  International Conference on}, pages 4144--4148. IEEE, 2016.

\bibitem[Lan and Zhou(2016)]{lan2016conditional}
G.~Lan and Y.~Zhou.
\newblock Conditional gradient sliding for convex optimization.
\newblock \emph{SIAM Journal on Optimization}, 26\penalty0 (2):\penalty0
  1379--1409, 2016.

\bibitem[Mokhtari et~al.(2018{\natexlab{a}})Mokhtari, Hassani, and
  Karbasi]{mokhtari2017conditional}
Aryan Mokhtari, Hamed Hassani, and Amin Karbasi.
\newblock Conditional gradient method for stochastic submodular maximization:
  Closing the gap.
\newblock In \emph{AISTATS}, pages 1886--1895, 2018{\natexlab{a}}.

\bibitem[Mokhtari et~al.(2018{\natexlab{b}})Mokhtari, Hassani, and
  Karbasi]{mokhtari2018stochastic}
Aryan Mokhtari, Hamed Hassani, and Amin Karbasi.
\newblock Stochastic conditional gradient methods: From convex minimization to
  submodular maximization.
\newblock \emph{arXiv preprint arXiv:1804.09554}, 2018{\natexlab{b}}.

\bibitem[Mokhtari et~al.(2018{\natexlab{c}})Mokhtari, Ozdaglar, and
  Jadbabaie]{mokhtari2018escaping}
Aryan Mokhtari, Asuman Ozdaglar, and Ali Jadbabaie.
\newblock Escaping saddle points in constrained optimization.
\newblock In \emph{Advances in Neural Information Processing Systems}, pages
  3629--3639, 2018{\natexlab{c}}.

\bibitem[Nguyen et~al.(2017{\natexlab{a}})Nguyen, Liu, Scheinberg, and
  Tak{\'a}{\v{c}}]{nguyen2017sarah}
Lam~M Nguyen, Jie Liu, Katya Scheinberg, and Martin Tak{\'a}{\v{c}}.
\newblock Sarah: A novel method for machine learning problems using stochastic
  recursive gradient.
\newblock In \emph{Proceedings of the 34th International Conference on Machine
  Learning-Volume 70}, pages 2613--2621. JMLR. org, 2017{\natexlab{a}}.

\bibitem[Nguyen et~al.(2017{\natexlab{b}})Nguyen, Liu, Scheinberg, and
  Tak{\'a}{\v{c}}]{nguyen2017stochastic}
Lam~M Nguyen, Jie Liu, Katya Scheinberg, and Martin Tak{\'a}{\v{c}}.
\newblock Stochastic recursive gradient algorithm for nonconvex optimization.
\newblock \emph{arXiv preprint arXiv:1705.07261}, 2017{\natexlab{b}}.

\bibitem[Reddi et~al.(2016)Reddi, Sra, P{\'o}czos, and
  Smola]{reddi2016stochastic}
Sashank~J Reddi, Suvrit Sra, Barnab{\'a}s P{\'o}czos, and Alex Smola.
\newblock Stochastic frank-wolfe methods for nonconvex optimization.
\newblock In \emph{2016 54th Annual Allerton Conference on Communication,
  Control, and Computing (Allerton)}, pages 1244--1251. IEEE, 2016.

\bibitem[Sahu et~al.(2018)Sahu, Zaheer, and Kar]{sahu2018towards}
Anit~Kumar Sahu, Manzil Zaheer, and Soummya Kar.
\newblock Towards gradient free and projection free stochastic optimization.
\newblock \emph{arXiv preprint arXiv:1810.03233}, 2018.

\bibitem[Shen et~al.(2019{\natexlab{a}})Shen, Fang, Zhao, Huang, and
  Qian]{shen2019complexities}
Zebang Shen, Cong Fang, Peilin Zhao, Junzhou Huang, and Hui Qian.
\newblock Complexities in projection-free stochastic non-convex minimization.
\newblock In \emph{The 22nd International Conference on Artificial Intelligence
  and Statistics}, pages 2868--2876, 2019{\natexlab{a}}.

\bibitem[Shen et~al.(2019{\natexlab{b}})Shen, Ribeiro, Hassani, Qian, and
  Mi]{shen2019hessian}
Zebang Shen, Alejandro Ribeiro, Hamed Hassani, Hui Qian, and Chao Mi.
\newblock Hessian aided policy gradient.
\newblock In \emph{International Conference on Machine Learning}, pages
  5729--5738, 2019{\natexlab{b}}.

\bibitem[Staib and Jegelka(2017)]{staib2017robust}
Matthew Staib and Stefanie Jegelka.
\newblock Robust budget allocation via continuous submodular functions.
\newblock In \emph{ICML}, pages 3230--3240, 2017.

\bibitem[Sutton and Barto(2018)]{sutton2018reinforcement}
Richard~S Sutton and Andrew~G Barto.
\newblock \emph{Reinforcement learning: An introduction}.
\newblock MIT press, 2018.

\bibitem[Vondr{\'a}k(2008)]{vondrak2008optimal}
Jan Vondr{\'a}k.
\newblock Optimal approximation for the submodular welfare problem in the value
  oracle model.
\newblock In \emph{STOC}, pages 67--74. ACM, 2008.

\bibitem[Wang et~al.(2016)Wang, Sadhanala, Dai, Neiswanger, Sra, and
  Xing]{wang2016parallel}
Yu-Xiang Wang, Veeranjaneyulu Sadhanala, Wei Dai, Willie Neiswanger, Suvrit
  Sra, and Eric Xing.
\newblock Parallel and distributed block-coordinate frank-wolfe algorithms.
\newblock In \emph{International Conference on Machine Learning}, pages
  1548--1557, 2016.

\bibitem[Yurtsever et~al.(2019)Yurtsever, Sra, and
  Cevher]{yurtsever2019conditional}
Alp Yurtsever, Suvrit Sra, and Volkan Cevher.
\newblock Conditional gradient methods via stochastic path-integrated
  differential estimator.
\newblock In \emph{International Conference on Machine Learning}, pages
  7282--7291, 2019.

\bibitem[Zhang et~al.(2019)Zhang, Chen, Mokhtari, Hassani, and
  Karbasi]{zhang2019quantized}
Mingrui Zhang, Lin Chen, Aryan Mokhtari, Hamed Hassani, and Amin Karbasi.
\newblock Quantized frank-wolfe: Communication-efficient distributed
  optimization.
\newblock \emph{arXiv preprint arXiv:1902.06332}, 2019.

\bibitem[Zheng et~al.(2018)Zheng, Bellet, and Gallinari]{zheng2018distributed}
Wenjie Zheng, Aur{\'e}lien Bellet, and Patrick Gallinari.
\newblock A distributed frank--wolfe framework for learning low-rank matrices
  with the trace norm.
\newblock \emph{Machine Learning}, 107\penalty0 (8-10):\penalty0 1457--1475,
  2018.

\bibitem[Zhou et~al.(2018)Zhou, Xu, and Gu]{zhou2018stochastic}
Dongruo Zhou, Pan Xu, and Quanquan Gu.
\newblock Stochastic nested variance reduced gradient descent for nonconvex
  optimization.
\newblock In \emph{Advances in Neural Information Processing Systems}, pages
  3921--3932, 2018.

\end{thebibliography}

\end{document}